\documentstyle{article}
\title{Notes on Seiberg--Witten Gauge Theory}
\author{Matilde Marcolli}
\date{June 1995; revised August 1995}
\newtheorem{teor}{Theorem}[section]
\newtheorem{defin}[teor]{Definition}
\newtheorem{lem}[teor]{Lemma}
\newtheorem{esem}[teor]{Example}
\newtheorem{corol}[teor]{Corollary}
\newtheorem{conj}[teor]{Conjecture}

\newcommand{\nc}{\newcommand}
\nc{\beq}{\begin{equation}}
\nc{\eeq}{\end{equation}}

\begin{document}

\maketitle

\section{Introduction}

$\Phi\upsilon\sigma\iota\varsigma$
$\kappa\rho\upsilon\pi\tau\varepsilon\sigma\vartheta\alpha\iota$
$\varphi\iota\lambda\varepsilon\iota$

(Eraclitos)

\vskip .2in

In the fall of 1994 E. Witten announced a ``new gauge theory of
4-manifolds'',  capable
of giving results analogous to the earlier theory of Donaldson (see
\cite{DK}), but where the computations involved are ``at least a thousand
times easier'' (Taubes).

The new theory begins with the introduction of the
monopole equations, whose physical motivation lies in the new results in $N=2$
supersymmetric Yang-Mills theory announced in \cite{WS1}, \cite{WS2}.

The equations are in terms of a section of a Spinor bundle and a
$U(1)$ connection on a line bundle $L$.
The first equation just says that the spinor section $\psi$ has to
be in the kernel of the Dirac operator. The second
equation describes a relation between the self-dual part of the
curvature associated to the connection $A$ and the section $\psi$ in
terms of the Clifford action.

The mathematical setting for Witten's gauge theory is considerably
simpler than Donaldson's analogue: first of all it deals with
$U(1)$-principal bundles (hermitian line bundles) rather than with
$SU(2)$-bundles, and the abelian structure group allows simpler
calculations; moreover the equation, which plays a role somehow
analogous to the previous anti-self-dual equation for
$SU(2)$-instantons (see \cite{DK}), involves Dirac operators and
$Spin_c$-structures, which are well known and long developed
mathematical tools (see \cite{R} or \cite{LM}).

The main differences between the two theories arise when it comes to
the properties of the moduli space of solutions of the monopole
equation up to gauge transformations. The
Seiberg--Witten invariant, which depends on the Chern class of the line
bundle $L$, is given by the number of points, counted with
orientation, in a zero-dimensional moduli space.

The present introduction to the subject of
Seiberg--Witten gauge theory will be taken mainly from the paper by
Witten \cite{W};
however extensive use will be made of other references that have
recently appeared.

In Witten's paper \cite{W} the monopole equation is introduced, and
the main properties of the moduli space of solutions are
deduced.

The dimension of the moduli space is computed by an index theory
technique, following an analogous proof for Donaldson's theory, as
in \cite{AHS};
and the circumstances under which the Seiberg--Witten invariants
provide a topological invariant of the four--manifold are illustrated
in a similar way to the analogous result regarding the Donaldson
polynomials.

The main conjecture, explained and justified in \cite{W} only at the
level of physical intuition, is that the Seiberg-Witten invariants
should coincide with invariants recently
derived by Kronheimer and Mrowka from the Donaldson polynomials
(\cite{KM2} and \cite{D}).

The tool that is of primary importance in proving the results about
the moduli space of abelian instantons is the Weitzenb\"{o}ck formula
for the Dirac operator on the $Spin_c$-bundle $S^+\otimes L$: such a
formula is a well known (see \cite{R}) decomposition of the square of
the Dirac operator on a spin bundle twisted with a line bundle $L$.

A first property which follows from the Weitzenb\"{o}ck formula is a
bound on the number of solutions: the moduli space is empty for all
but finitely many choices of the line bundle $L$.

Moreover, as shown in \cite{KM1}, the moduli space is always compact:
a fact that avoids the complicated analytic techniques that were
needed for the compactification of the moduli space of
$SU(2)$-instantons (see \cite{DK}).

Another advantage of this theory is that the singularities of the
moduli space (again this is shown in
\cite{KM1}) only appear at the trivial section $\psi\equiv 0$, since
elsewhere the action of the gauge group is free. Hence,
by perturbing the equation, it is possible to get a smooth moduli space.

The analogue of the vanishing theorem for
Donaldson polynomials on a manifold that splits as a connected sum
can be proven, reinforcing the intuitive feeling that the two sets of
invariants ought to be the same.

Moreover explicit computations can be done in the case of K\"{a}hler
manifolds, by looking at the
Seiberg-Witten invariants associated to the canonical line bundle.

The latter result has a generalization due to Taubes \cite{T}, where
it is shown that the value $\pm1$ of the Seiberg-Witten invariants is
achieved on symplectic four--manifolds, with respect to the canonical
line bundle, by a technique that involves
estimates of solutions of a parametrized family of perturbed monopole
equations.

The meaning of such a result (as Taubes recently explained in a series
of lectures held at the AMS meeting in San Francisco) lies in the fact
that, once compared with the vanishing theorem, it seems to
indicate that symplectic manifolds are somehow the most ``basic''
(indecomposable) kind of four-manifolds: a point of view that leads
to new insights on the geometry of four-manifolds.

It is still to be understood what the good conjecture should be; in fact,
as proven recently in \cite{K}, the strongest conjecture formulated by
Taubes is not true.

Another important application of the new gauge theory has been the
positive answer to the Thom conjecture for embedded surfaces in ${\bf
  C}P^2$, proved in \cite{KM1}: an oriented two manifold $\Sigma$ embedded in
${\bf C}P^2$ and representing the same homology class as an algebraic
curve of degree $d$ has genus $g$
\[ g\geq \frac{(d-1)(d-2)}{2}. \]

These notes do not intend to be a comprehensive
introduction to the Seiberg--Witten gauge theory: something of that
sort will be available, I believe, in a very short time, and written
by someone who is a real expert in the field. What's written here may
be better thought of as a guide on ``how to approach the study of
Seiberg--Witten gauge theory''. It collects the notes prepared for a
series of seminars given at the University of Milano, and intended
mainly for an audience of first, or second year graduate students, or advanced
undergraduate students in Mathematics or Physics who prepare their
dissertation in topology, differential geometry, or quantum field theory.
The purpose of these notes is to help the student in reading the
references avilable on the subject.
As prerequisites we shall only assume some knowledge of the topology
of fibre bundles, and some notions like that of connection, and curvature,
that are usually covered in an advanced undergraduate course of algebraic
topology, or differential geometry; as well as of theoretical, or
mathematical physics. Any other concept will be briefly introduced,
and the reader who is not familiar with the subject will be addressed
to some references, where she/he can find more detailed information.
In particular, Seiberg--Witten gauge theory is considered here as a
completely self--contained subject: no a priori knowledge of Donaldson
theory is assumed, and any mention of that can be skipped by the
reader who is not familiar with the non--abelian gauge theory.
We warn the reader that the last chapters, starting with the quantum
field theory, are rather tentative: in fact not much has been written
up yet. Therefore our introduction will be very sketchy,
and in part just based on speculation; however, we hope it will help
to give a feeling of how the current research is oriented.

A revised edition of these notes is in preparation. It will contain
supplementary chapters on Seiberg--Witten--Floer homology, and on
the Gromov Invariants and applications to Quantum Cohomology.
Some other chapters will be rewritten with more details, and more
up to date material.

\noindent{\bf Acknowledgement.}
I am deeply grateful to Mel Rothenberg, for having helpfully and patiently
supervised my study of gauge theory; I also wish to thank
Clifford Taubes for having sent me some of his preprints and notes. I
greatly enjoyed the course on Donaldson's gauge theory
given by Simon Salamon at the Summer School of Mathematics in Cortona and the
lectures on the Seiberg--Witten invariants given by Taubes at
the AMS Conference in San Francisco. I thank J.Peter May for his continuous
support and encouragement. Last but not least, I thank
Renzo Piccinini, for having invited me to give these seminars at the
University of Milano; and all those who attended, for their
intelligent comments and questions, and for their appreciation.

\vspace{.5cm}

\section{Preliminary Notions}

We shall introduce here some of the basic concepts that are
needed in order to define the Seiberg--Witten equations and
invariants; this introduction will be rather sketchy: more detailed
information can be found in any of the standard references listed in
the bibliography, such as \cite{LM} or \cite{R}.

\subsection{Clifford Algebras and Dirac Operators}

\begin{defin}
The Clifford algebra $C(V)$ of a (real or complex)  vector space $V$ with a
symmetric  bilinear
form $(,)$ is the algebra generated by the elements
\[ e_1^{\epsilon_1}\cdots e_n^{\epsilon_n}, \]
where $\epsilon_i=0$, or 1 and $\{ e_i \}$ is an orthogonal basis of
$V$, subject to the relations
\[ e\cdot e^\prime + e^\prime \cdot e =-2(e,e^\prime). \]
The multiplication of elements of $V$ in the Clifford algebra is
called Clifford multiplication.
\end{defin}

In particular given a differentiable manifold $X$ we shall consider
the Clifford algebra associated to the tangent space at each point.

\begin{defin}
The Clifford algebra of the tangent bundle of $X$ is the bundle that has
fibre over each point $x\in X$ the Clifford algebra $C(T_x X)$. We
shall denote this bundle $C(TX)$.
\end{defin}

If $dim V= 2m$, there is a unique irreducible representation of the
Clifford algebra $C(V)$. This representation has
dimension $2^m$.

\begin{defin}
A Spinor Bundle over a Riemannian manifold  $X$ is
the vector bundle associated to $C(TX)$ via this irreducible representation,
endowed with a hermitian structure such that the Clifford multiplication is
skew--symmetric and  compatible with the Levi--Civita connection on
$X$ (see \cite{R}).
\label{spinors}
\end{defin}

Not all manifolds admit a spinor bundle; it is proved in \cite{R} that
the existence of such a bundle is equivalent to the existence of a
$Spin_c$--structure on the manifold $X$: we shall discuss
$Spin_c$--structures in the next paragraph.
If such a bundle exists, it splits as a direct sum of two
vector bundles,
\beq
S = S^+ \oplus S^-,
\label{Splus}
\eeq
where the splitting is given by the internal grading of the Clifford
algebra.

As a consequence of definition \ref{spinors} we can prove an easy and
useful lemma.

\begin{lem}
Let $X$ be a manifold that admits a Spinor bundle $S$. Let $\{ e_i \}$ be
a local orthonormal basis of sections of the tangent bundle $TX$
and $<,>$ be the
hermitian structure as in definition \ref{spinors}. Then, for any
section $\psi\in \Gamma(X,S)$, the expression
\[ <e_ie_j\psi, \psi > \]
is purely imaginary at each point $x\in X$.
\label{real}
\end{lem}

\noindent\underline{Proof:} In fact, by skew--adjointness of Clifford
multiplication and the fact that the basis is orthonormal,
\[ <e_ie_j\psi,\psi >=-<e_j\psi,e_i\psi >=<\psi,e_je_i\psi>
-<\psi,e_ie_j\psi>=-\overline{<e_ie_j\psi,\psi >}. \]

\noindent QED

\begin{defin}
Given a spinor bundle $S$ over $X$, the Dirac operator on $S$ is a
first order  differential operator on the smooth sections
\[ D: \Gamma (X, S^+)\rightarrow \Gamma (X, S^-) \]
defined  as the composition
\beq
D: \Gamma (X, S^+) \stackrel{\nabla}{\rightarrow} \Gamma
(X, S^+) \otimes T^*X
\stackrel{g}{\rightarrow}  \Gamma (X, S^+)  \otimes
TX  \stackrel{\bullet}{\rightarrow}  \Gamma (X, S^-),
\label{D}
\eeq
where the first map is the covariant derivative, with the
$Spin$--connection induced by the Levi--Civita
connection on $X$ (see \cite{R}), the second is the
Legendre transform given by the Riemannian metric, and the third is
Clifford multiplication.
\end{defin}

It is easy to check (for more details see \cite{R}) that this
corresponds to the following expression in coordinates:
\[ Ds=\sum_k e_k \cdot \nabla_k s. \]

An essential tool in Spin geometry, which is very useful in
Seiberg--Witten gauge theory (see e.g. \cite{JPW},
\cite{KM1},\cite{T}, \cite{W}), is the {\em Weitzenb\"ock formula}.

\begin{defin}
Given a smooth vector bundle $E$ over a $Spin_c$--manifold $X$ and a
connection $A$ on $E$,
the twisted Dirac operator $D_A: \Gamma (X, S^+\otimes E) \rightarrow
\Gamma (X, S^-\otimes E)$ is the operator acting on a section
$s\otimes e$ as the Dirac operator on $s$ and the composite of the
covariant derivative $\tilde\nabla_A$ and the Clifford multiplication on $e$:
\[ D_A: \Gamma(X, S^+\otimes E)\stackrel{\nabla\otimes
  1+1\otimes\tilde\nabla_A}{\rightarrow} \Gamma(X, S^+\otimes E\otimes
T^*X) \]
\[ \stackrel{g}{\rightarrow} \Gamma(X, S^+\otimes E\otimes TX)
\stackrel{\bullet}{\rightarrow} \Gamma(X, S^-\otimes E).  \]
\label{twist}
\end{defin}

Note that in the above, as will be explained in the following section,
the bundle $S^+\otimes E$ may exist
even in cases when neither of the two bundles $S^\pm$ and $E$ exist
separately.

The twisted Dirac operator is still
well defined, as it is defined in terms of local quantities;
and the following theorem holds true as
well (see \cite{R}, \cite{LM}).

\begin{teor}
The twisted Dirac operator $D_A$ satisfies the
Weitzenb\"ock formula:
\[ D^2_A s = (\nabla_A^*\nabla_A + \frac{\kappa}{4} + \frac{-i}{4}F_A)s, \]
where $\nabla_A^*$ is the formal adjoint of the covariant derivative
with respect to the $Spin$--connection on the Spinor bundle, and with
respect to the connection $A$ on $E$, $\nabla_A=\nabla\otimes 1 + 1 \otimes
\tilde\nabla_A$; $\kappa$ is the scalar curvature on $X$,
$F_A$ is the curvature of the connection $A$, and $s\in \Gamma (X,
S^+\otimes E)$ .
\label{weitzenbock}
\end{teor}

\noindent\underline{Proof:}In local coordinates
\[ D_A^2 s=\sum_{ij} e_i {\nabla_A}_i(e_j {\nabla_A}_j s)=
\sum_{ij}e_i e_j {\nabla_A}_i {\nabla_A}_j s = \]
\[ = -\sum_i {\nabla_A}_i^2 s + \sum_{i<j} e_i e_j
({\nabla_A}_i {\nabla_A}_j-{\nabla_A}_j {\nabla_A}_i) s, \]
where $\nabla_A=\nabla\otimes 1 +1\otimes \tilde\nabla_A$. The first summand
is $\nabla_A^*\nabla_A$ (\cite{R} pg.28), and the second splits into a term
which corresponds to the scalar curvature on $X$ (\cite{R} pg. 126) and
the curvature $-iF_A$ of the connection $A$ (we identify the
Lie algebra of $U(1)$ with $i{\bf R}$).

\noindent QED

\subsection{$Spin$ and $Spin_c$ Structures}

The group $Spin(n)$ is the universal covering of $SO(n)$.

The group $Spin_c(n)$ is defined via the following extension:
\beq
1\rightarrow {\bf Z}_2\rightarrow Spin_c(n)\rightarrow SO(n)\times
U(1) \rightarrow 1,
\label{Spinc}
\eeq
i.e. $Spin_c(n)=(Spin(n)\times U(1))/{\bf Z}_2$.

The extension (\ref{Spinc}) determines the exact
sheaf--cohomology sequence:

\beq
\cdots \rightarrow H^1(X;Spin_c(n))\rightarrow H^1(X;SO(n))\oplus
H^1(X;U(1))\stackrel{\delta}{\rightarrow} H^2(X;{\bf
  Z}_2)
\label{cohom}
\eeq

By the standard fact that $H^1(X;G)$ represents the equivalence
classes of principal $G$--bundles over $X$, we see that the connecting
homomorphism  of the sequence (\ref{cohom}) is given by
\[ \delta : (P_{SO(n)},P_{U(1)})\mapsto w_2(P_{SO(n)})+\bar
c_1(P_{U(1)}), \]
where $\bar c_1(P_{U(1)})$ is the reduction mod 2 of the first Chern
class of the line bundle associated to the principal bundle $P_{U(1)}$
by the standard representation and $w_2$ is the second
Stiefel--Whitney class.

\begin{defin}
A manifold $X$ has a $Spin_c$--structure if the frame bundle lifts
to a principal $Spin_c(n)$ bundle. It has a $Spin$--structure if it
lifts to a $Spin(n)$ principal bundle.
\label{Spinstr}
\end{defin}

{}From the above considerations on the cohomology sequence
(\ref{cohom}), and analogous considerations on the group $Spin(n)$,
we can prove a useful lemma.

\begin{lem}
A manifold $X$ admits a $Spin_c$--structure iff $w_2(X)$ is the
reduction mod 2 of an integral class. It has a $Spin$--structure iff
$w_2(X)=0$.
\end{lem}

\begin{corol}
Different $Spin_c$--structures on $X$ are parametrised by
\[ 2H^2(X;{\bf Z})\oplus H^1(X;{\bf Z}_2). \]
\end{corol}

This follows directly form (\ref{cohom}): see \cite{LM}.

With this characterisation of $Spin_c$--structures we have the
following theorem, which is proved in \cite{HH}.

\begin{teor}
Every oriented 4--manifold admits a $Spin_c$--structure.
\end{teor}

\subsection{Spinor Bundles}

Both $Spin(n)$ and $Spin_c(n)$ can be thought of as lying inside the
Clifford algebra
$C({\bf R}^n)$, \cite{R}, \cite{LM}. Therefore to a principal
$Spin(n)$ or $Spin_c(n)$ bundle we can associate a vector bundle via
the unique irreducible representation of the Clifford algebra.
This will be the bundle of Spinors over $X$
associated to the $Spin_c$ or $Spin$ structure, as defined in
\ref{spinors}.

This can be given an explicit description in terms of transition
functions (see \cite{LM}). In fact let $g_{\alpha\beta}$ be the
transition functions of the frame bundle over $X$, which take values in
$SO(n)$. Then locally they can be lifted to functions $\tilde
g_{\alpha\beta}$ which take values in $Spin(n)$, since on a differentiable
manifold it is always possible to choose open sets with contractible
intersections that trivialise the bundle.

However, if we have a $Spin_c$--manifold that is not $Spin$, the
$\tilde g_{\alpha\beta}$ won't form a cocycle, as
$ \tilde g_{\alpha\beta} \tilde g_{\beta\gamma} \tilde
g_{\gamma\alpha} =1$
means exactly that the second Stiefel--Whitney class vanishes.

Because of the $Spin_c$ structure we know that $w_2$ is the reduction
of an integral class $c\in H^2(X,{\bf Z})$, which represents
a complex line bundle, say with transition functions
$\lambda_{\alpha\beta}$ with values in $U(1)$. Such functions will
have a square root $\lambda_{\alpha\beta}^{1/2}$ locally;
however, the line bundle will not have a square root globally (which
is to say that the $\lambda_{\alpha\beta}^{1/2}$ won't form a
cocycle), since by construction the first Chern class is not divisible by
2.

However, the relation $w_2(X)+c=0$ mod 2 that comes from (\ref{cohom})
says that the product
\beq
\tilde g_{\alpha\beta}   \lambda_{\alpha\beta}^{1/2}
\label{trans}
\eeq
is a cocycle. These are the transition functions of $S\otimes
L$, where $S$ would be the Spinor bundle of a $Spin$ structure and $L$
would be the square root of a line bundle: neither of these objects is
defined globally, but the tensor product is.
This is the description of the Spinor Bundle of a $Spin_c$ structure
that we shall use in the following.

\subsection{Symplectic and K\"{a}hler Manifolds}

Recall that a manifold $X$ is endowed with a symplectic structure if a
closed 2--form $\omega$ is given on
$X$ such that its highest exterior power is nowhere vanishing.

The symplectic structure determines a compatible almost complex structure $J_x:
T_xX\rightarrow T_xX$, $J_x^2=-1$, and inner product, by representing
the skew--symmetric bilinear form $\omega_x(v,w)=<v, J_x w>_x$.
When this distribution $J_x$ is
integrable, in the sense of the Frobenius theorem, we have an actual
complex structure on $X$. In that case $X$ is a K\"{a}hler manifold.
In particular $g(v,w)=\omega(v,Jw)$ is a Riemannian metric on $X$
which is compatible with the complex structure.

It is known that all symplectic structures are locally the same. In
this sense symplectic geometry can be thought of as something more
rigid than ${\cal C}^\infty$--geometry but less rigid than Riemannian
geometry. In fact by the Darboux theorem (see \cite{Fo}) it is always
possible to find a set of coordinates in which the symplectic form
reduces to the ``standard one'': $\omega=\sum_i dx_i\wedge dy_i$ (a
symplectic manifold must be even dimensional).

A typical example of symplectic manifold is the cotangent bundle of a
smooth manifold $Y$, $X=T^*Y$ with the local coordinates $\{ q_i, p_i
\}$ and the symplectic form $\omega=\sum_i dp_i\wedge dq_i$: this is
the mathematical setup of classical Hamiltonian Mechanics.

Note that, unlike symplectic geometry, which is less rigid than Riemannian
geometry, K\"ahler geometry is much more rigid, which is expected in
passing from smooth to analytic geometry. In this sense, the condition
of being K\"ahler is rather exceptional, although it may be non-trivial
to provide examples of compact symplectic manifolds that are not
K\"ahler, particularly in the simply connected case; see
e.g. \cite{McD}.

A useful notion is that of {\em symplectic reduction}.
The quotient of a symplectic manifold with respect to
the action of a group is not in general a symplectic manifold; however
under some suitable hypothesis it is possible to describe the orbit space
as a union of symplectic manifolds. A good reference for this topic
is \cite{By}.

The action of a group on a symplectic manifold is said to be symplectic if it
preserves the form $\omega$.

\begin{defin}
A symplectic action of a Lie group $G$ on $(X,\omega)$ is Hamiltonian if
each vector field $v$ on $X$, spanned by the infinitesimal action of the Lie
algebra ${\cal L}\rightarrow Vect(X)$, lifts to a map
$H_v\in {\cal C}^\infty (X)$, via the relation
\[ \omega(v,\cdot)=dH_v (\cdot). \]
\label{hamilt}
\end{defin}

\begin{defin}
The {\em moment map} of a Hamiltonian symplectic action of $G$ on $X$ is
the map
\[ \mu : X \rightarrow {\cal L}^* \]
\[  x \mapsto (\mu (x) : {\cal L} \rightarrow {\bf R}) \]
\[  \mu (x) (v) = H_v (x). \]
\label{moment}
\end{defin}

\begin{teor}
Suppose given a Hamiltonian symplectic action of a group $G$ on $(X,\omega)$
and an element
$\xi \in {\cal L}^*$ such that $\mu^{-1}(\xi)\subset X$ is a submanifold.
Then $X_\xi =\mu^{-1}(\xi)/ G_\xi$ is a manifold, where $G_\xi$
is the stabiliser
of $\xi$ under the induced coadjoint action; moreover $X_\xi$ has a symplectic
structure induced via the pullback to $\mu^{-1}(\xi)$ of the form $\omega$.
\label{reduction}
\end{teor}

Symplectic geometry plays a prominent role in the Seiberg--Witten gauge
theory. Many
computations are possible for the case of symplectic 4--manifolds
(\cite{T}, \cite{T1}, \cite{KMT});
moreover, although the
invariants are defined in terms of $U(1)$--connections, and sections
of Spinor bundles, as we shall illustrate below, they turn out to be
strictly related to invariants of symplectic manifolds, known as
Gromov invariants (see \cite{McDS}, \cite{T2}).
The fact that the
Seiberg--Witten invariants have a ``more basic'' structure also led to
a conjecture, suggested by Taubes, that symplectic manifolds may be
among the most basic building blocks of the whole geometry of
4--manifolds.

\subsection{The Index Theorem}

We shall recall here very briefly some essential results of Index Theory.
The reader who is not familiar with the topological and analytic
properties of the index of elliptic operators is urged to gain some
familiarity with the Atiyah--Singer Index Theorem. A good and
very readable source of information is the book by R. Boos and
D.D. Bleecker, \cite{BB}.

A differential operator of order $m$, mapping the smooth sections of
a vector bundle $E$ over a compact manifold $Y$ to those of another
such bundle $F$, can be described in local
coordinates and local trivialisations of the bundles as
\[ D=\sum_{\mid\alpha\mid\leq m} a_\alpha(x)D^\alpha, \]
with $\alpha=(\alpha_1,\ldots,\alpha_n)$. The coefficients $a_\alpha(x)$ are
matrices of smooth functions that represent elements of $Hom(E,F)$ locally;
and $D^\alpha=\frac{\partial}{\partial x^{\alpha_1}_1}\cdots
\frac{\partial}{\partial x^{\alpha_n}_n}$.

\begin{defin}
The principal symbol associated to the operator $D$ is the expression
\[ \sigma_m(D)(x,p)=\sum_{\mid\alpha\mid=m}a_\alpha(x)p^\alpha. \]
\label{sym}
\end{defin}

Given the differential operator $D:\Gamma(Y,E)\rightarrow\Gamma(Y,F)$,
the principal symbol with the local expression above defines a
global map
\[ \sigma_m : \pi^*(E)\rightarrow \pi^*(F), \]
where $T^*Y\stackrel{\pi}{\to}Y$ is the cotangent bundle; that is, the
variables $(x,p)$ are local coordinates on $T^*Y$.

Consider bundles $E_i$, $i=1\ldots k$, over a compact $n$-dimensional
manifold $Y$ such
that there is a complex $\Gamma(Y,E)$ formed by the spaces of
sections $\Gamma(Y,E_i)$
and differential operators $d_i$ of order $m$:
\[ 0\rightarrow \Gamma(Y,E_1)\stackrel{d_1}{\rightarrow}\cdots
\stackrel{d_{k-1}}{\rightarrow}\Gamma(Y,E_k)\rightarrow 0. \]

Construct the principal symbols $\sigma_m(d_i)$; these determine an associated
symbol complex
\[ 0\rightarrow \pi^*(E_1)\stackrel{\sigma_m(d_1)}{\rightarrow}\cdots
\stackrel{\sigma_m(d_{k-1})}{\rightarrow}\pi^*(E_k)\rightarrow 0.\]

\begin{defin}
The complex $\Gamma(Y,E)$ is elliptic iff the associated symbol complex
is exact.
\label{ellip}
\end{defin}

In the case of just one operator, this means that $\sigma_m(d)$ is an
isomorphism off the zero section.

The Hodge theorem states that the cohomology of the complex $\Gamma(Y,E)$
coincides with the harmonic forms; that is,
\[ H^i(E)=\frac{Ker(d_i)}{Im(d_{i-1})}\cong Ker(\Delta_i), \]
where $\Delta_i=d^*_id_i+d_{i-1}d^*_{i-1}$.

Without loss of generality, by passing to the assembled complex
\[ E^+=E^1\oplus E^3\oplus\cdots \]
\[ E^-=E^2\oplus E^4\oplus\cdots, \]
we can always think of one elliptic operator $D:\Gamma(Y,E^+)
\rightarrow\Gamma(Y,E^-)$, $D=\sum_i (d_{2i-1}+d^*_{2i})$.

The Index Theorem is as follows:
\begin{teor}
Consider an elliptic complex over a compact, orientable, even dimensional
manifold $Y$ without boundary. The index of $D$, which is given by
\[ Ind(D)=\dim Ker(D)-\dim Coker(D)=\sum_i (-1)^i\dim Ker \Delta_i
=-\chi(E), \]
$\chi(E)$ being the Euler characteristic of the complex, can be expressed
in terms of characteristic classes as:
\[ Ind(D)=(-1)^{n/2}<\frac{ch(\sum_i (-1)^i [E_i])}{e(Y)}td(TY_{\bf C}),
[Y]>. \]
\label{indthm}
\end{teor}

In the above $ch$ is the Chern character, $e$ is the Euler class of the
tangent bundle of $Y$, $td(TY_{\bf C})$ is the Todd class of the complexified
tangent bundle.

To give more details on the Index Theorem would imply introducting $K$--theory
with compact support, pseudodifferential operators, and many other tools.
We leave it to the reader to discover all the beauty of the Atiyah--Singer
theorem, by reading through the appropriate references.

\subsection{Other useful notions}

Recall that the gauge group of a $G$--bundle is defined as  the group
of self equivalences of the
bundle, namely the group of smooth maps
\[ \lambda_\alpha : U_\alpha \rightarrow G \]
\[ \lambda_\beta = g_{\beta\alpha}\lambda_\alpha g_{\alpha\beta}, \]
where the bundle is trivial over $U_\alpha$ and has transition
functions $g_{\alpha\beta}$.
It is often useful to consider this space of maps endowed with Sobolev
norms, and, by completing with respect to the norm, to consider gauge
groups of $L^2_k$ functions, as in \cite{DK}, or \cite{FU}.

The gauge group is an infinite dimensional
manifold; if the structure group is abelian then it has a simpler
description as ${\cal G}={\cal M}(X, G)$, the space of maps from $X$
to $G$.

\begin{lem}
In the case $G=U(1)$, the set of connected components of the
gauge group ${\cal G}$ is $H^1(X,{\bf Z})$.
\label{conngauge}
\end{lem}

\section{The Functional and the Equations}

In all the following we shall consider $X$ to be a compact, connected,
orientable,  differentiable 4--manifold without boundary. We shall
refer to a $Spin_c$--structure on $X$ by specifying the Spinor Bundle
$S\otimes L$.

\subsection{The Equations}

The equations of the gauge theory are given in terms of a pair $(A,\psi)$ of
indeterminates, of which $A$ is a connection on $L$ and
$\psi$ is a smooth section of $S^+\otimes L$.

The equations are
\beq
D_A\psi=0
\label{eqSW1}
\eeq
\beq
(F^+_A)_{ij}=\frac{1}{4}<e_ie_j\psi,\psi>e^i\wedge e^j,
\label{eqSW2}
\eeq
where $D_A$ is the Dirac operator twisted by the connection $A$, and
$F^+_A$ is the self--dual part of the curvature associated to
$A$. Here $\{ e_i \}$ is a local basis of $TX$ that acts on $\psi$ by
Clifford multiplication (see the exercises), $\{ e^i \}$ is the dual
basis of $T^*X$, and $<,>$ is the inner product on the fibres of
$S^+\otimes L$.

We follow the notation of \cite{JPW} rather than that of \cite{W}
or \cite{KM1}; the equivalence of these notations is left as an
exercise at the end of this chapter.

\subsection{The Gauge Group}

The gauge group of $L$ is well defined although $L$ is not globally
defined as a line bundle, since the definition of the gauge group is given
just in terms of the transition functions. In particular, as in the
case of a line bundle, ${\cal G}={\cal M}(X, U(1))$.

There's an action of the gauge group on the space of pairs $(A,\psi)$,
where $A$ is a
connection on $L$ and $\psi$ a section of $S^+\otimes L$, given by
\beq
\lambda: (A, \psi) \mapsto (A-2i \lambda^{-1}d\lambda, i\lambda \psi).
\label{action}
\eeq

\begin{lem}
The action defined in (\ref{action}) induces an action of ${\cal G}$
on the space of solutions to the Seiberg--Witten equations.
\end{lem}

\noindent\underline{Proof:} It's enough to check that
\[ D_{A-2i \lambda^{-1}d\lambda}(\lambda \psi)=i\lambda D_A\psi +
id\lambda \cdot \psi -2id\lambda \cdot \psi +id\lambda\cdot\psi, \]
and in the second equation
\[ F^+_{A-2i \lambda^{-1}d\lambda}=F^+_A -2d^+(i
\lambda^{-1}d\lambda)= F^+_A, \]
and $<e_ie_j\lambda \psi,\lambda\psi>=\mid \lambda\mid^2 <e_ie_j\psi,\psi>=
<e_i e_j\psi,\psi>$.

\noindent QED

It is clear from (\ref{action}) that the action of ${\cal G}$ on the
space of solutions is free iff $\psi$ is not identically zero; while for
$\psi\equiv 0$ the stabiliser of the action is $U(1)$, the group of
constant gauge transformations.

\subsection{The Seiberg--Witten Functional and the Variational Problem}

Given a moduli problem formulated in terms of differential equations,
one may consider a functional of which the equations represent the
absolute minima: an example is the case of the Yang--Mills functional,
and the anti--self--dual equation for $SU(2)$ Donaldson gauge theory.

In the case of the Seiberg--Witten equations, it's not hard to figure
out what such a functional could be: in fact it suffices to consider
the following.

\begin{defin}
The Seiberg--Witten functional acting on pairs $(A,\psi)$ is given by
\beq
S(A,\psi)=\int_X (\mid D_A\psi \mid^2 + \mid
F^+_A-\frac{1}{4}<e_ie_j\psi,\psi>e^i\wedge e^j \mid^2) dv.
\label{functional}
\eeq
\end{defin}

\begin{lem}
Via the Weitzenb\"ock formula (theorem \ref{weitzenbock}) the
Seiberg--Witten functional (\ref{functional}) can be rewritten as
\beq
S(A,\psi)=\int_X (\mid \nabla_A\psi \mid^2 + \mid F_A^+ \mid^2 +
\frac{\kappa}{4} \mid \psi \mid^2 + \frac{1}{8}\mid \psi \mid^4) dv.
\label{functional2}
\eeq
\end{lem}

\noindent\underline{Proof:}
In fact we have that
\[ \mid D_A\psi \mid^2= <D_A^2\psi, \psi >= <\nabla^*_A\nabla_A\psi,
\psi > + \frac{\kappa}{4}\mid\psi\mid^2 + \frac{1}{4}<F_A\psi, \psi > \]
and
\[  \mid F^+_A-\frac{1}{4}<e_ie_j\psi,\psi>e^i\wedge e^j \mid^2 = \]
\[ = \mid F_A^+ \mid^2 -\frac{1}{4}(F_A^+,<e_ie_j\psi,\psi> e^i\wedge e^j)+
\frac{1}{16}\mid <e_ie_j\psi,\psi>\mid^2 \mid  e^i\wedge e^j\mid^2, \]
where $(,)$ denotes the pointwise inner product of two--forms:
$(\alpha,\beta)dv= \alpha\wedge *\beta$. But $(F_A^+, e^i\wedge
e^j)={F_A^+}_{ij}$ and $\frac{1}{4}(F_A^+,<e_ie_j\psi,\psi>
e^i\wedge e^j)=<\frac{1}{4}{F_A^+}_{ij}e_ie_j\psi,\psi>$, which is the
expression of the action of $F_A^+$ on $\Gamma(X,S^+\otimes L)$, via
Clifford multiplication.
Thus two terms cancel out in the sum, as in the first
summand of (\ref{functional}) only the self--dual part of the
curvature acts non--trivially on the section $\psi$ (see the exercise
at the end of the section). Moreover $\mid
e^i\wedge e^j\mid^2=1$ and $\mid <e_ie_j\psi,\psi>\mid^2=2\mid
\psi\mid^4$.

\noindent QED

It seems clear from the definition that solutions of (\ref{eqSW1}) and
(\ref{eqSW2}) are the absolute minima of (\ref{functional}). Note however,
that, in order to make this claim precise, we would need the property of
coercivity of the functional, which is shown later in this section.

Some properties that follow from introducing the Seiberg--Witten
functional are summarised in the following lemma.

\begin{lem}
If the scalar curvature of $X$ is non-negative, all solutions of the
Seiberg--Witten equations have $\psi\equiv 0$. Under the assumption
that $X$ has $b^{2^+}>1$ for a generic choice of
the metric, the only solutions will be $\psi=0$ and $A$ flat.
\label{poscurv}
\end{lem}

\noindent\underline{Proof:} Clearly (\ref{functional}) and the
Weitzenb\"ock formula imply $\psi=0$. They also give $F^+_A=0$; but
the first Chern class of the line bundle $L$, which is given by
$c_1(L)=\frac{1}{2\pi}F_A$, is an integral class modulo torsion; and
this implies that $F_A\in H^{2-}(X;{\bf R})\cap H^2(X;{\bf Z})/T$,
with $T$ the torsion subgroup of $H^2(X;{\bf Z})$. But $H^2(X;{\bf Z})/T$ is a
lattice in $H^2(X;{\bf R})$, and $H^{2-}(X;{\bf R})$ is a subspace of
codimension strictly greater than 1 (as $b^{2^+}>1$); therefore for a
generic choice of the metric, which means that this subspace
$H^{2-}(X;{\bf R})$ is in generic position with respect to the
lattice, they will not intersect outside the point $F=0$. Thus the
connection $A$ is a flat connection.

\noindent QED

The above argument fails in the case $b^{2^+}=1$, as in that case, even
for generic metrics, there would be non--trivial intersections of the
lattice of integral forms with the anti-self-dual harmonic forms.

Clearly the following holds true as well.

\begin{corol}
If, moreover, the first Chern class of $L$ in $H^2(X;{\bf Z})$ is not a
torsion element, the above lemma gives that the only solution is the
trivial one $A=0$, $\psi=0$.
\end{corol}

In the case where $c_1(L)$ is torsion, we can still get rid of the
flat connections by suitably perturbing the equations.

\begin{defin}
The perturbed Seiberg--Witten equations are obtained by introducing a
small harmonic self--dual two--form $\eta$ as a perturbation parameter:
\[  D_A\psi=0, \]
\[  F_A^+ +i\eta = \frac{1}{4}<e_ie_j\psi,\psi >e^i\wedge e^j. \]
\label{perturb}
\end{defin}

The corresponding perturbed functional will be
\[ S_\eta (A,\psi)=\int_X (\mid F^+_A \mid^2 + \mid \nabla_A\psi \mid^2
+ \frac{\kappa}{4}\mid \psi \mid^2)dv + \]
\[ \int_X F^+_A\wedge i\eta +\int_X
\mid \frac{1}{4}<e_ie_j\psi,\psi>e^i\wedge e^j -i\eta\mid^2 dv. \]

It is clear from the perturbed functional that the corresponding
equations, given in definition \ref{perturb}, cannot have solutions
with flat connections, as that would correspond to
\[ \mid \frac{1}{4}<e_ie_j\psi,\psi>e^i\wedge e^j-i\eta\mid =0, \]
which is not compatible with $\psi=0$ and $\eta\neq 0$. Therefore we
have the following.

\begin{corol}
For a manifold $X$ with $b^{2+} >1$ the perturbed
Seiberg--Witten equations will not have solutions with $\psi=0$.
\end{corol}

In particular this means that the corresponding moduli space will not
have singular points.

It's a general fact in gauge theory that, given a functional like
(\ref{functional}), one may look at the absolute minima, or just at
the extremals, i.e. at solutions of the Euler--Lagrange equations.
The equation for minima is in general a first order problem, while the
Euler--Lagrange equations will give a second order problem.
For instance, the functional considered in Donaldson theory is the
Yang--Mills functional. The anti--self--dual connections are
the absolute minima; while the corresponding variational problem
gives rise to the Yang--Mills equation (see \cite{T5}).

In introducing the variational problem for the Seiberg--Witten
functional (\ref{functional}), and the analytic properties of
(\ref{functional}), we follow \cite{JPW}.

\begin{lem}
The Euler--Lagrange equations are of the form
\beq
D_A^2\psi-\frac{i}{2}F_A^+\cdot
\psi-\frac{1}{8}<e_ie_j\psi,\psi>e_ie_j\psi=0
\label{var1}
\eeq
and
\beq
d^*(F^+_A-\frac{1}{4}<e_ie_j\psi,\psi>e^i\wedge
e^j)+\frac{i}{2}Im<D_A\psi, e_i\psi>e^i=0.
\label{var2}
\eeq
\end{lem}

Again by using the Weitzenb\"ock formula we can rewrite these
equations as

\[ \nabla_A^*\nabla_A\psi
+\frac{\kappa}{4}\psi+\frac{1}{4}\mid\psi\mid^2\psi =0 \]
and
\[ d^*F_A^+ +\frac{i}{2}Im<\nabla_i\psi,\psi>e^i=0. \]

An interesting question (as already pointed out in
\cite{JPW}) is what the geometric meaning of the other
non--minimising solutions could be.

\subsection{Analytic properties of the Seiberg--Witten functional}

Here we introduce some analytic properties of the Seiberg--Witten
functional, following \cite{JPW}. These will be the main tool in
proving the compactness of the moduli space, which is the major
difference between Seiberg--Witten and Donaldson gauge theory.

The first property that we shall introduce is that the functional
(\ref{functional}) is coercive; hence we'll see that it has the
Palais--Smale condition. We leave it to an upcoming revised version of
these notes to deal with other more geometric properties, as the
Mountain--pass  lemma and the conditions under which the
Seiberg--Witten functional is a Morse--Bott function.

The latter property is particularly interesting in the case of the
Seiberg--Witten functional of
the reduced theory on three manifolds (see the chapter on dimensional
reduction): in fact it is
the first step needed in order to define a homology, via the Morse theory
of Seiberg--Witten monopoles. At the moment Taubes is working on this
Seiberg--Witten--Floer homology \cite{T3}
while there is already some work done from the Physicist's point of view
\cite{BL}.

\vspace{.3cm}

In our context the analytic data will be $L^2_1(S^+\otimes L)$ and
$L^2_1({\cal A})$,
where the first denotes the completion in the $L^2_1$--norm of the space of
smooth sections $\Gamma (X,S^+\otimes L)$ and the second is the
$L^2_1$--completion of
the space of 1--forms associated to the affine space of connections, given a
fixed connection $A_0$. We shall also consider, as gauge transformations,
the space ${\cal G}^2_2$, which is locally the completion of ${\cal G}$
in the $L^2_2$--norm. These are all infinite dimensional Hilbert manifolds
\cite{JPW}.

There's an action of ${\cal G}^2_2$ on $L^2_1(S^+\otimes L)\times
L^2_1({\cal A})$,
which is differentiable and is given by (\ref{action}).

Consider the Seiberg--Witten functional on the space
$L^2_1(S^+\otimes L)\times L^2_1({\cal A})$. It is proven in \cite{JPW}
that this functional is coercive, i.e. the following holds.

\begin{lem}
There is a constant $c$ such that, for some $g\in {\cal G}^2_2$,
\[ S(A,\psi)\geq c^{-1}(\| \lambda\psi \|_{L^2_1} +\| A-2i\lambda^{-1}
d\lambda \|_{L^2_1})-c. \]
\label{coerc}
\end{lem}

\subsubsection{The Palais--Smale Condition}

Recall the following definition \cite{JPW}:

\begin{defin}
Suppose given a Hilbert manifold $Y$ with the action of a (possibly
infinite dimensional) Lie group $G$ and a smooth functional
$f:Y\rightarrow {\bf R}$,
which is $G$ invariant. Then $f$ has the Palais--Smale condition if,
for any sequence $\{ x_k \}\subset Y$ such that

(i) $f(x_k)$ is bounded and

(ii) $df(x_k)\to 0$, as $k\to\infty$,

\noindent there is a subsequence $\{ x^\prime_k \}$ and a sequence
$\{ g_k \}\subset G$ such that $g_kx^\prime_k\to x$, with $df(x)=0$ and
$f(x)=\lim_k f(x_k)$.
\label{PS}
\end{defin}

In order to prove that the Seiberg--Witten functional over the space
$L^2_1(S^+\otimes L)\times L^2_1({\cal A})$ has the Palais--Smale condition,
we need some preliminary results.

\begin{lem}
Let $(A,\psi)$ be a solution of the Seiberg--Witten equations.
Then the following estimate holds:
\[ \|\psi\|_{L^\infty}^2\leq \max (0,-\min_{x\in X}\kappa(x)). \]
\label{psibound}
\end{lem}

\noindent\underline{Proof:}Note that, since we are considering the equations
in $L^2_1$--spaces, we are not assuming any regularity.

We want to show that the set $S$ of points where $\mid\psi\mid >1$ is of
measure
zero.

Use the Weitzenb\"ock formula together with the Euler--Lagrange equation
(\ref{var1}) to get
\[ \nabla^*_A\nabla_A \psi = -(\frac{\kappa}{4}\mid\psi\mid^2-\frac{1}{8}
\mid\psi\mid^4). \]
Introduce a function $\eta$ such that
$\eta=(\mid\psi\mid -1)\frac{\psi}{\mid\psi\mid}$ if $\mid\psi\mid >1$,
and $\eta=0$ if $\mid\psi\mid \leq 1$.
A straightforward computation (we are following the argument given in
\cite{JPW}) gives that
\[ 0\geq\int_S <\nabla_A\psi,\nabla_A\eta >+ \frac{1}{8}(\mid\psi\mid^2-1)
(\mid\psi\mid -1)\mid\psi\mid , \]
assuming that the scalar curvature $\kappa \geq -1$, and that
\[ <\nabla_A\psi,\nabla_A\eta>\geq 0; \]
hence $S$ is of measure zero.

\noindent QED

A somehow simpler proof, which requires however some a priori assumption of
regularity of the solutions, can be found in \cite{KM1}.

\begin{lem}
Fix a connection $A_0$. Given a sequence of connections $\{ A_k \} \subset
L^2_1({\cal A})$, there exists a sequence of gauge transformations
$\{ \lambda_k \}$ in the identity component of ${\cal G}^2_2$
such that the 1--forms
\[ A_k -2i\lambda_k d\lambda_k -A_0 \]
are co--closed.
\label{coclosed}
\end{lem}

\noindent\underline{Proof:}This can be read off from the complex
(\ref{chcomplex}), which will be introduced in the next chapter.
In fact, the directions in the tangent space to the connections that
satisfy the Seiberg--Witten equations, which are spanned by the action of
the gauge group, are in the image of the exterior derivative $d$ (because
of the expression of the infinitesimal action of ${\cal G}$). Hence, by
acting on $A_k$ with a gauge transformation $\lambda_k$ which is connected to
the identity of the gauge group, we can obtain that
$A_k -2i\lambda_k d\lambda_k$ is in the orthogonal complement, with respect
to the chosen norm, to the image of $d$, i.e. in the kernel of $d^*$.

\noindent  QED

\begin{teor}
Given a sequence $\{ (A_k,\psi_k) \}$ of solutions of the Seiberg--Witten
equations in
$L^2_1(S^+\otimes L)\times L^2_1({\cal A})$, there exist a subsequence
$\{ (A_{k'},\psi_{k'}) \}$
and a sequence of gauge transformations $\{ \lambda_{k'} \}$
in ${\cal G}^2_2$, such that the sequence
\[ \{ (A_{k^\prime}-2i\lambda_{k^\prime} d\lambda_{k^\prime},
 i\lambda_{k^\prime}\psi_{k^\prime}) \} \]
converges with all derivatives to a solution $(A,\psi)$ of the
Seiberg--Witten equations.
\label{comp1}
\end{teor}

\noindent\underline{Proof:}Using lemma \ref{coclosed} we can assume that
$\{ A_k -A_0 \}$ is a sequence of co--closed 1--forms.
Moreover, since $F_{A_k}=dA_k$, we have $dF_{A_k}=0$.
Thus, $d^* F_{A_k}= (d+d^*) F_{A_k}$
and $d^*dA_k=(d^*d + dd^*)A_k$, where both $d+d^*$ and $d^*d+dd^*$ are
elliptic operators.

Hypothesis (i) of definition \ref{PS} holds, since
$S(A_k,\psi_k)\equiv 0$, and, since the
Seiberg--Witten functional is coercive, this means that
$\| \lambda_k\psi \|_{L^2_1}$, and
$\| A_k-2i\lambda_k d\lambda_k-A_0 \|_{L^2_1}$
are bounded (possibly after composing with another family of gauge
transformations).

The bound on the $L^2_1$--norm of the connections gives an
$L^2$ bound:
\[ \| A_k-A_0 \|_{L^2}\leq c \| A_k- A_0 \|_{L^2_1}. \]

In order to show that there is a subsequence that converges with all
derivatives, it is sufficient to show that all Sobolev norms are bounded
and use the Sobolev embedding theorem.

To bound higher Sobolev norms observe that
\[ d^* F_{A_k} = 2 d^* F_{A_k}^+ = -i\sum_j Im(<\nabla_j\psi_k,\psi_k>)e^j, \]
from the variational equation (\ref{var2}), and use the facts that
$*d* F=*dF^+ -*dF^-$ and $0=*dF=*dF^+ +*dF^-$.
This gives a bound on $\| d^* F_A \|_{L^2}$ and therefore on
$\| d^*d A_k \|_{L^2}$.

Now the elliptic estimate applied to the operator $dd^*+d^*d$ says that
\[ \| A_k \|_{L^2_2}\leq c( \| d^*d A_k \|_{L^2} + \| A_k \|_{L^2}). \]
Its is clear that this procedure can be carried over for all higher
Sobolev norms.

The result for the sections $\psi_k$ follows from lemma \ref{psibound}, which
gives an $L^2$ bound, and the coercivity property of the Seiberg--Witten
functional, which gives the $L^2_1$ bound. The bounds on the higher
norms are obtained by applying the elliptic estimate to the Dirac operator.

\noindent QED

The result of theorem \ref{comp1} holds true also for solutions of the
variational problem (\ref{var1}), (\ref{var2}). Moreover, the argument given
here can be slightly modified in order to show the following result
(whose proof is given in \cite{JPW}).

\begin{teor}
The Seiberg--Witten functional (\ref{functional}) satisfies the Palais--Smale
condition of definition \ref{PS}.
\label{PS2}
\end{teor}

\subsection{Exercises}
\begin{itemize}

\item The equation (\ref{eqSW2}) is often written in
  terms of endomorphisms of $S^+$.
  Following \cite{DK} we can consider the map $\gamma : TX \rightarrow
  Hom(S^+,S^-)$ given in local coordinates by $\gamma (e_i)=\sigma_i$, where
  the $\sigma_i$ are the Pauli
  matrices; $S^+$ and $S^-$  are complex two-plane bundles.
  This induces an action of the two--forms $\Lambda^2$ on $S^+$ given by
  the expression
  $e\wedge e' (s)= -\gamma^*(e)\gamma(e')s$,
  where $*$ denotes the adjoint (\cite{DK} pg. 76).
  Check that $\Lambda^{2-}$ acts trivially,
  and therefore this can be considered as an action of the self--dual
  two--forms. Hence we get a map
  $\rho : \Lambda^{2+}\rightarrow End(S^+)$.
  Check that (\ref{eqSW2}) can be written as
  \[ \rho(F^+_A)=\sigma (\psi\otimes \bar\psi), \]
  where $\sigma$ is the projection on the traceless part of
  $End(S^+\otimes L)\cong S^+\otimes S^+$.

\item Check that the map $\rho$ in the above problem changes norms by
  a factor 2.

\item The equations (\ref{eqSW1}), (\ref{eqSW2}) can be written as
\[ D_A \psi =0 \]
and
\[ (F^+_A)_{ij}= -\frac{i}{2}\bar\psi\Gamma_{ij}\psi, \]
where the $\Gamma_{ij}=[\sigma_i,\sigma_j]$ are defined in terms of the
Clifford matrices (see \cite{W}).
Since $ \psi\in \Gamma (X, S^+\otimes L)$,
$\bar\psi \in \Gamma (X, S^+\otimes L^{-1})$ and thus the right
hand side can also be read as
\[ \bar\psi\cdot\psi \in \Gamma (X, S^+\otimes S^+)
\cong Hom(S^+,S^+). \]
Check that this sheaf splits as
\[ Hom(S^+,S^+)\cong \Lambda^0(TX)\oplus \Lambda^{2+}(TX) \]
due to the Clifford action, and that, in the second equation,
the self-dual part of the curvature has to coincide with the
projection of $\bar\psi\cdot \psi$ on
$\Lambda^{2+}(TX)$.

\item Complete the proof of lemma \ref{psibound}.

\item Prove that $\mid <e_ie_j\psi,\psi>\mid^2=2\mid\psi\mid^4$.

\end{itemize}

\section{Seiberg--Witten Invariants}

\subsection{The Moduli Space}

We can consider the moduli space of solutions of (\ref{eqSW1}),
(\ref{eqSW2}), modulo the action of the gauge group. The purpose of
studying the topology of the moduli space is to have a somehow ``simpler''
model of the manifold $X$ by means of which to compute invariants
associated to the differentiable structure of $X$.

\begin{defin}
The Moduli Space $M$ is the set of solutions of the Seiberg--Witten
equations, modulo the action of the gauge group.
$M$ depends on the choice of the
$Spin_c$ structure, i.e. of the line bundle $L^2$.
\label{moduli}
\end{defin}

In principle, by the above definition, it seems that we are considering
a whole, possibly infinite, family of moduli spaces, according to the
choice of $L^2$; however in the following we'll see (as proved in
\cite{W}) that only finitely many choices will give rise to a
nontrivial moduli space.

Some of the most important geometric properties of the moduli space
are finite dimensionality, compactness, orientability, and the fact
that it has at most very ``nice'' singularities.
We are going to give a proof of these results in the rest of this
section, following mainly \cite{KM1} and \cite{W}.

\subsubsection{Computation of the Dimension}

One of the most striking features of the Atiyah--Singer Index Theorem
is certainly the fact that it provides a powerful tool for the sort of
computations as the one that we shall illustrate in this section. In
fact, in many different contexts (see
\cite{AHS},\cite{DK},\cite{McDS},\cite{W})
where one wishes to compute the dimension of the moduli space of
solutions of certain differential equations modulo the action of some
large symmetry group, one tries to construct a local model of the moduli space
by linearising the equations to some Fredholm operator and then fit
the linearisation into a short chain complex such that its Euler
characteristic, computed via the Index Theorem, gives the dimension of
the moduli space.

This argument works in some generic case (i.e for an open dense set of
metrics); therefore, for
metrics that are ``non--generic'' in that sense, the moduli space may
have a dimension that is not the expected one. This is the case,
for Seiberg--Witten equations,
when one considers K\"{a}hler metrics.

Thus, in order to compute the dimension of $M$, we shall linearize
(\ref{eqSW1}) and (\ref{eqSW2}) in a neighbourhood of a pair
$(A_0,\psi_0)$.

\begin{lem}
\label{linear}
The linearized Seiberg--Witten equations at a pair
$(A_0+i\alpha,\psi_0+\phi)$, with $\alpha$ a 1--form, and $\phi$ a
section of $S^+\otimes L$,  are
\[ D_{A_0}\phi + i\alpha\cdot \psi_0 =0 \]
and
\[ d^+\alpha -\frac{1}{2}Im(<e_ie_j\psi_0,\phi>)e^i\wedge e^j=0. \]
\end{lem}

Now consider the infinitesimal action of the gauge group on the
solution $(A_0,\psi_0)$. If we write an element of the
gauge group as a map $\lambda=e^{i f}$ for some $f:X\rightarrow
{\bf R}$, this means that we consider
\[ (A_0,\psi_0)\mapsto (A_0 -idf, if\psi_0), \]
by our definition of the action of ${\cal G}$, (\ref{action}).

Now consider the following short sequence of spaces and maps:
\beq
0\rightarrow\Lambda^0 \stackrel{G}{\rightarrow} \Lambda^1\oplus
\Gamma(S^+\otimes L) \stackrel{T}{\rightarrow}
\Lambda^{2+}\oplus\Gamma(S^+\otimes L)\rightarrow 0.
\label{chcomplex}
\eeq
Here $G$ is the map given by the infinitesimal action of ${\cal G}$;
$T$ is the operator defined by the linearization of the
Seiberg--Witten equations, i.e. the left hand side of the equations in
lemma \ref{linear}; the $\Lambda^q$ are $q$--forms on $X$ and
$\Gamma(S^+\otimes L)$ is the space of smooth sections of the Spinor Bundle.

\begin{lem}
The sequence (\ref{chcomplex}) is a chain complex. We shall indicate
this complex by $C^*$. The operators $T$ and $G$ are Fredholm.
\end{lem}

\noindent\underline{Proof:} We need to check that $T\circ G=0$. But in
fact
\[ D_{A_0}(if\psi_0) -idf\cdot \psi_0 =0, \]
because of (\ref{eqSW1}), and
\[ d^+(df)+\frac{1}{2}Im(<e_ie_j\psi_0,if\psi_0>)e^i\wedge e^j = \]
\[ =Im(\frac{if}{2}<e_ie_j\psi_0,\psi_0>)e^i\wedge e^j =0. \]
Here we used the facts that $d^+d=p_{\Lambda^{2+}}\circ d^2=0$ and that
$<e_ie_j\psi_0,\psi_0>$ is purely imaginary (as proved in lemma \ref{real}).

The operators $T$ and $G$ are Fredholm since, up to zero--order terms,
they are given by the elliptic differential operators
$d^+ +d^*$ and $D_A$.

\noindent QED

By definition the tangent space of $M$ at the point $(A_0,\psi_0)$ is
the quotient $Ker(T)/Im(G)$: in fact we consider the linear approximation
to the Seiberg--Witten equations modulo those directions that are
spanned by the action of the gauge group.

Thus we need to compute $H^1(C^*)$. The Index Theorem provides a way
to compute the Euler characteristic of $C^*$ in terms of some
characteristic classes.

\begin{teor}
The Euler characteristic of $C^*$ is
\[ -\chi(C^*)= Ind (D_A+d^+ +d^*), \]
where $d^*$ is the adjoint of the exterior derivative.
\label{eulchar}
\end{teor}

\noindent\underline{Proof:} Up to zero--order operators
$G$ can be deformed to the exterior derivative $d$;
and $T$ can be deformed to the pair of operators
\[ D_A: \Gamma (X, S^+\otimes L) \rightarrow\Gamma (X, S^-\otimes L) \]
and
\[ d^+: \Lambda^1\rightarrow \Lambda^{2+}. \]
Hence the assembled complex becomes
\[ 0\rightarrow\Lambda^1\oplus \Gamma (X, S^+\otimes L)
\stackrel{D_A+d^+ +d^*}{\rightarrow}
\Lambda^0\oplus \Lambda^{2+}\oplus \Gamma (X, S^-\otimes L)
\rightarrow 0. \]
By the Index Theorem the Euler Characteristic of the original complex
is not affected by this change, hence:
\[ -\chi (C^*)=Ind(D_A+d^+ +d^*). \]

\noindent QED

\begin{corol}
The index of the complex $C^*$ is equal to
\[ c_1(L)^2-\frac{2\chi+3\sigma}{4}, \]
where $\chi$ is the Euler characteristic of $X$, $\sigma$ is the
signature of $X$, and $c_1(L)^2$ is the cup product with itself of the
first Chern class of $L$ integrated over $X$ (with a standard abuse of
notation we write $c_1(L)^2$ instead of $<c_1(L)^2, [X]>$).
\end{corol}

\noindent\underline{Proof:} Use the additivity of the index. By the
index theorem for the twisted Dirac operator (\cite{R}, \cite{BB}) it
is known that
\[ Ind(D_A)=-\int_X ch(L)\hat A(X). \]
The Chern character is $ch(L)=2(1+c_1(L)+\frac{1}{2}c_1(L)^2+\cdots)$
(the rank of $L$ over the reals
is two). The $\hat A$ class is $\hat A(X)=1-\frac{1}{24}p_1(X)+\cdots$,
where $p_1(X)$ is the first Pontrjagin class of the tangent bundle.
Thus the top
degree term of $ch(L)\hat A(X)$ will be $\frac{1}{12}p_1(X)+c_1(L)^2$.
On the other hand, by the index theorem for the signature operator
(again see \cite{R}, or \cite{BB}) it is known that $\frac{1}{3}\int_X
p_1(X)= \sigma$, hence we get
\[ Ind(D_A)=c_1(L)^2-\frac{\sigma}{4}. \]

Note that, since $L$ is not really a line bundle, its first Chern
class is defined to be $c_1(L)=\frac{c_1(L^2)}{2}$, and it makes
sense in the coefficient ring ${\bf Z}[\frac{1}{2}]$.

The index of $d^+ +d^*$ can be read off from the chain complex
\[ 0\rightarrow \Lambda^0 \stackrel{d}{\rightarrow} \Lambda^1
\stackrel{d^+}{\rightarrow} \Lambda^{2+}\rightarrow 0. \]
The Euler characteristic of this complex turns out to be
\[ -Ind(d^*+d^+)=\frac{1}{2}(\chi+\sigma), \]
by another application of the Index Theorem.

Summing up together it follows that
\[ \chi (C^*)=c_1(L)^2-\frac{2\chi+3\sigma}{4}. \]

\noindent{QED}

Now, in order to obtain from this index computation the dimension of
the moduli space, we need the following lemmata:

\begin{lem}
The tangent space $T_{(A,\psi)}M$ at a regular point $\psi\neq 0$ can
be identified with $H^1(C^*)$.
\label{tgsp}
\end{lem}

\noindent\underline{Proof:} $H^1(C^*)$ describes exactly those
directions that
are spanned infinitesimally at the point $(A,\psi)$ by the solutions
of the Seiberg--Witten equations, modulo those directions that are
spanned by the action of the gauge group.

\noindent QED

\begin{lem}
In the above complex, $H^0(C^*)=0$, and, under a suitable perturbation
of the Seiberg--Witten equations, also $H^2(C^*)=0$.
\end{lem}

\noindent\underline{Proof:} $H^0(C^*)=0$ since the map $G$, which describes
the infinitesimal action of the gauge group, as in (\ref{chcomplex}), is
injective. In order to show that $H^2(C^*)=0$, consider the perturbed
Seiberg--Witten equations, as in definition \ref{perturb}, and the
corresponding deformation
\[ 0\rightarrow \Lambda^0 \stackrel{G}{\rightarrow}
\Lambda^1\oplus\Lambda^{2+}\oplus
\Gamma(S^+\otimes L) \stackrel{\tilde T}{\rightarrow}
\Lambda^{2+}\oplus\Gamma(S^+\otimes L)\rightarrow 0 \]
of (\ref{chcomplex}). The map $\tilde T$ is given as
\[ \tilde T (\alpha,\phi,\eta)=(D_{A_0}\phi +
i\alpha\cdot \psi_0, d^+\alpha +\eta
-\frac{1}{2}Im<e_ie_j\psi_0,\phi>e^i\wedge e^j). \]

Suppose there is a section in $\Lambda^{2+}\oplus\Gamma(S^+\otimes
L)$, that is orthogonal to the image of $\tilde T$ in the space of
$L^2$ sections. We want to show that such an element must be
identically zero. By an argument of \cite{R}, this is enough to
prove that $\tilde T$ is surjective on smooth sections as well.
So assume that $(\chi,\beta)$ is orthogonal to any section of the form
\[ (D_{A_0}\phi + i\alpha\cdot \psi_0, d^+\alpha +\eta
-\frac{1}{2}Im<e_ie_j\psi_0,\phi>e^i\wedge e^j). \]
This means that
\[ <D_{A_0}\phi + i\alpha\cdot \psi_0,\chi>=0, \]
with the inner product of sections of $S^-\otimes L$, and
\[ (d^+\alpha +\eta -\frac{1}{2}Im<e_ie_j\psi_0,\phi>e^i\wedge e^j,
\beta)=0 \]
as two--forms. By the arbitrariness of $\alpha$ and $\eta$, both
$\chi$ and $\beta$ must be zero.

A cobordism argument \cite{DK} shows that any such perturbation of the
equations does not affect the invariants.

\noindent QED

The above lemmata and computations yield the following result.
\begin{teor}
The dimension of the moduli space $M$ at a regular point is given by
\[ dim(M)= c_1(L)^2-\frac{2\chi+3\sigma}{4}. \]
\label{dim}
\end{teor}

The result proven here is an infinitesimal result; the technique for
passing to a local result is the same used in the $SU(2)$ gauge
theory: as a reference we can point to \cite{AHS}.

\subsubsection{Compactness}

Everything needed in order to show that the moduli space is compact
is contained in theorem \ref{comp1}. In fact we proved there that
every sequence of solutions of the Seiberg--Witten equations has a
convergent subsequence, up to gauge transformations. This implies the
following

\begin{teor}
Any sequence of points in the moduli space of connections and
sections $(A,\psi)$ that satisfies the Seiberg--Witten equations, modulo
the action of the gauge group ${\cal G}$, has a convergent subsequence.
Therefore, the moduli space is compact.
\label{comp2}
\end{teor}

\subsubsection{Orientation}

The proof of the orientability of the moduli space given in \cite{W}
just mimics an analogous argument in Donaldson theory, \cite{D2}: the
orientation is given by a trivialisation of the determinant line
bundle associated to the linearized Seiberg--Witten equations.

There's another proof of orientability given for Donaldson theory,
\cite{FU}, which doesn't go through to the Seiberg--Witten case
unless some stronger assumptions are made on the manifold $X$, for example,
that $X$ be simply connected. Though this proof is not as good, it's
a sufficient simplification for our purposes,
therefore we shall present it briefly.

Following \cite{FU}, note first that the critical points of $M$
do not affect the orientability: therefore consider the non-compact
manifold
\[ \hat M=M-\{ singularities \}. \]
There's an embedding of $\hat M$ in the space
$\hat{\cal A}$ of pairs $(A,\psi)$ with $\psi\neq 0$:
\[ \hat M \hookrightarrow \hat{\cal A}/{\cal G}. \]
The action of ${\cal G}$ on $\hat{\cal A}$ is free.

It is known (see e.g. \cite{BB}) that a family of Fredholm
operators on a space $Y$ determines an {\em index bundle} $[Ind(D)]\in
KO(Y)$, in
$K$--theoretic language. by definition $[ Ind(D) ]$ is orientable iff
the characteristic class $w_1 ([Ind(D)])=0$. Recall that, for an element
in $KO(Y)$, $w_1(\xi - \zeta):= w_1(\xi)-w_1(\zeta)$ in ${\bf Z}_2$.

Moreover, if the space $Y$ is simply connected, then every bundle over $Y$
is orientable.

Thus, the strategy to prove the orientability of the moduli space is to
realize the tangent bundle $T\hat M$ as a subbundle of the index bundle
of a family of Fredholm operators over a simply connected space.

According to lemma \ref{tgsp}, the fibre of $T\hat M$ over $(A,\psi)$ is
given by $H^1(C^*)$. In other words, the complex $C^*$ corresponds to
an assignment of a Fredholm operator to the point $(A,\psi)$.
Hence we get a family of Fredholm operators parametrised by all possible
choices of $(A,\psi)$, $\psi\neq 0$, namely the linearization of the
Seiberg--Witten equations at the chosen connection and section modulo
the action of the gauge group.
By theorem \ref{eulchar} this family of Fredholm operators may be thought
of as the following:
\[ T_{(A,\psi)}=D_A + d^+ + d^* . \]

To make this argument precise we should consider the appropriate Sobolev
norms on ${\cal G}$ and on the space of sections where the family $T$ acts.
We address the reader to \cite{FU}, where the analysis for $SU(2)$
gauge theory is developed.

Consider the index bundle of the above family $T$. By construction
$T\hat M$ is the pullback of $[Ind(T)]$ via the embedding of $\hat M$ in
$\hat{\cal A}/{\cal G}$.

Thus, we have the following:

\begin{teor}
Suppose that the manifold $X$ has $H^1(X;{\bf Z})=0$.
Then the moduli space $M$ is orientable.
\label{orient}
\end{teor}

\noindent\underline{Proof:}We want to show that the bundle
$T\hat M$ is orientable. By the above argument, it is enough to
show that $\hat{\cal A}/{\cal G}$ is simply connected.

Consider the fibration
\[ {\cal G}\rightarrow \hat{\cal A}\rightarrow \hat{\cal A}/{\cal G}. \]
The space of all pairs of a connection $A$ and a section $\psi$ is
contractible as product of an affine and a vector space. However, $\hat
{\cal A}$ might not in general be contractible since, imposing the
condition that
$\psi\neq 0$, we remove the origin of the vector space. But there is enough
codimension to have $\pi_1(\hat{\cal A})=0$.

Thus, in the long homotopy sequence of the fibration, we have
\[ \cdots \rightarrow \pi_1({\cal G})\rightarrow \pi_1(\hat{\cal A})
\rightarrow \pi_1(\hat{\cal A}/{\cal G})\rightarrow \pi_0({\cal G})
\rightarrow \pi_0(\hat{\cal A}); \]
where we have that $\pi_1(\hat{\cal A})=0$. Moreover, by lemma
\ref{conngauge}, $\pi_0{\cal G}=H^1(X;{\bf Z})$, which is trivial by
hypothesis. Hence $\pi_1(\hat{\cal A}/{\cal G})=0$.

\noindent QED

Note that, as already pointed out in the beginning of this paragraph, a
more general proof can be given of the orientability, which doesn't assume
any hypothesis on the cohomology of $X$. Note moreover that a technical
difficulty that we omitted to mention in the above argument arises in
considering the $K$--ring $KO(\hat{\cal A}/{\cal G})$, since
$\hat {\cal A}/{\cal G}$ is not a
compact space. We address the reader to \cite{FU}, or \cite{BB}, for
a more detailed treatment of these technicalities.

\subsection{The Invariants}

We shall assume here that $X$ has $b^{2+}>1$, to have the generic condition.

When the dimension
\[ c_1(L)^2 -\frac{2\chi +3\sigma}{4} < 0, \]
there are generically no solutions to the Seiberg--Witten equations.
Hence we define the Seiberg--Witten invariant $N_L$ to be zero.

If $L^2$ is such that the dimension is zero,
\[ c_1(L)^2 =\frac{2\chi +3\sigma}{4}, \]
generically the moduli space $M$ consists of a finite number of points, due
to compactness.
Since we also have an orientation of the moduli space, we can
associate to $M$ a number, which is obtained by counting the point
with a sign given by the orientation.

\begin{defin}
The Seiberg--Witten invariant, relative to a choice of $L$ such that
\[ c_1(L)^2-\frac{2\chi+3\sigma}{4}=0, \]
is given by
\[ N_L\equiv \sum_{p\in M} \epsilon_p, \]
with $\epsilon_p=\pm 1$ according to the orientation of $M$ at the
point $p$.
\label{inv0}
\end{defin}

It is more difficult to define the invariant when the dimension of the
moduli space is positive. In fact a priori there are many possible
ways to get a number by evaluating some cohomology class over the
cycle $M$, so as to generalise the counting of points.

We shall adopt the following definition (see \cite{T2}), which turns
out to be a good candidate with respect to some of the applications related
to symplectic geometry and quantum cohomology, that will be mentioned later
on.

In the section regarding the quantum field theoretic approach to
Seiberg--Witten theory we'll spend some more words on how to define
the invariants.

Consider the group of all gauge transformations that fix a base point,
${\cal G}_0\subset {\cal G}$. Take the moduli space $M_0$ of solutions of
Seiberg--Witten equations, modulo the action of ${\cal G}_0$.

\begin{lem}
This space $M_0$ fibres as a principal $U(1)$ bundle
over the moduli space $M$.
\end{lem}

Let ${\cal L}$ denote the line bundle over $M$ associated to this
principal $U(1)$ bundle via the standard representation. Then we
introduce the following invariant.

\begin{defin}
The Seiberg--Witten invariant, relative to a choice of $L$ such that
the dimension of $M$ is positive and even,
\[ d=c_1(L)^2-\frac{2\chi+3\sigma}{4}>0, \]
is given by the pairing of the $d/2$th-power of the Chern class of the
line bundle ${\cal L}$ with the moduli space $M$,
\[ N_L\equiv \int_M c_1({\cal L})^{d/2}. \]
If the dimension of $M$ is odd, the invariant is set to be zero.
\label{invpos}
\end{defin}

We'll show in the following section that only finitely many choices of
$L$ determine a nontrivial invariant.

\subsection{Finiteness}

Again this is going to be a consequence of the Weitzenb\"ock formula,
proved in theorem \ref{weitzenbock}. Consider the Seiberg--Witten functional
(\ref{functional}), rewritten in the form (\ref{functional2}).

\begin{lem}
Solutions to Seiberg--Witten equations have a uniform bound on $\int_X
\mid F^+_A \mid^2 dv$.
\end{lem}

\noindent\underline{Proof:}Complete
the term $\frac{\kappa}{4}\mid\psi\mid^2+\frac{1}{8}\mid\psi\mid^4$
to a square: this gives the estimate
\[ 0\leq \int_X (\frac{\kappa}{4}\mid\psi\mid^2+\frac{1}{8}\mid\psi\mid^4+
\frac{1}{8}\kappa^2)dv. \]

At a solution we have $S(A,\psi)=0$, hence
\[ \int_X \mid F^+_A \mid^2 dv =-\int_X (\mid \nabla_A\psi \mid^2 +
+ \frac{\kappa}{4}\mid\psi\mid^2+\frac{1}{8}\mid\psi\mid^4) dv. \]

Thus we get that
\[ \int_X \mid F^+_A \mid^2 dv \leq \int_X (\mid F^+_A \mid^2 +
\mid \nabla_A\psi \mid^2) dv = \]
\[ =-\int_X (\frac{\kappa}{4}\mid\psi\mid^2
+\frac{1}{8}\mid\psi\mid^4)dv \leq \frac{1}{8}\int_X \kappa^2 dv. \]

\noindent QED

Moreover, from the dimensional computation and the fact that
\[ c_1(L)^2=\int_X c_1(L)\wedge
c_1(L)=\frac{1}{(2\pi)^2}\int_X (\mid F^+_A\mid^2 -\mid F^-_A\mid^2)dv, \]
we also get the following.

\begin{lem}
Solutions to Seiberg--Witten equations have both $I^+=\int_X \mid
F^+_A \mid^2 dv$, and $I^-=\int_X \mid F^-_A \mid^2 dv$ uniformly
bounded, i.e. bounded by geometric quantities that do not depend on
the line bundle $L^2$.
\end{lem}

But there may be only finitely many choices of $L^2$ such that both
$I^+$ and $I^-$ are bounded (\cite{DK}, \cite{FU}, \cite{W}); hence
we have a main theorem of Seiberg--Witten gauge theory:
\begin{teor}
Only finitely many choices of the $Spin_c$ structure give rise to
non--trivial invariants.
\end{teor}

\subsection{A Cobordism Argument}

The reason we are considering the Seiberg--Witten invariants is
to introduce diffeomorphism invariants of a 4-manifold. Hence we need
to show that the construction above leads to a value of $N_L$ that is
independent of the metric on $X$. This turns out to be the case, at least
under the assumption that $b^{2+}(X)>1$. In the case $b^{2+}(X)=1$
the space of metrics breaks into {\em chambers}: inside each chamber the
choice of the metric doesn't affect the invariant; while when a path
of metrics crosses a {\em wall} between two chambers, the invariant
jumps by a certain amount. An accurate study of the structure of
chambers in the K\"ahler case can be found in \cite{FM}.

Here we shall present briefly an argument that shows the invariance
with respect to the metric for manifolds with $b^{2+}>1$. The proof of
the following results carries over to the present case from
\cite{DK} chapter 4 and 9. We shall
consider only the case when the moduli space is zero-dimensional, for
simplicity.

\begin{defin}
Consider a trivial infinite dimensional bundle over $\hat{\cal A}/{\cal
G}\times
{\cal C}$, where ${\cal C}$ is the space of metrics, with fibre
$\Lambda^{2+}\oplus \Gamma(X,S^-\otimes L)$.
Let $\sigma$ be a section (see lemma \ref{quant1} below)  given by
\[ \sigma (A,\psi)=(F^+_A-\frac{1}{4}<e_ie_j\psi,\psi>e^i\wedge e^j,
D_A\psi). \]
For a fixed $Spin_c$ structure on $X$,
we define the {\em universal moduli space} to be the zero set of $\sigma$:
\[ {\cal M}_L=\sigma^{-1}(0). \]
\label{univmod}
\end{defin}

In order to apply Donaldson's argument, it is essential to verify that
the following result holds in our case.

\begin{lem}
The projection map $\pi:{\cal M}_L\to{\cal C}$ is Fredholm.
\end{lem}

The next result allows us to generalise the
transversality theorem to an infinite dimensional context, \cite{DK}, pg.145.

\begin{lem}
If $f:{\cal X}\rightarrow {\cal Y}$ is a Fredholm map between infinite
dimensional Banach manifolds and $h: W\rightarrow {\cal Y}$ is a smooth map
from a finite dimensional manifold $W$, then $h$ can be arbitrarily
approximated
in the ${\cal C}^\infty$ topology by a map which is transverse to $f$.
\label{transverse}
\end{lem}

\begin{teor}
The invariants $N_L$ of a manifold $X$ with $b^{2+}(X)>1$ are
diffeomorphism invariants, independent of the metric.
\end{teor}

Suppose given two metrics $g_0$ and $g_1$, and a path $\gamma_t$
between them.
By the above, since $b^{2+}(X)>1$ we can choose $\gamma$ in such a way as to
avoid metrics that have singular solutions with $\psi\equiv 0$.

We want to show that the algebraic sum of the points in the zero set
$M_{L,g_t}$ of the section $\sigma(\cdot,\gamma_t)$ does not depend on $t$.

Consider the maps $\pi:{\cal M}_L\to {\cal C}$ and $\gamma: [0,1]\to {\cal C}$.
By the transversality theorem, lemma \ref{transverse}, the set
\[ M_L(\gamma_t)=\{ (A,\psi) \mid \pi (A,\psi,g)=g=\gamma_t \} \]
is a smooth 1-dimensional manifold with boundary.

{}From theorem \ref{comp2}, we deduce that $M_L(\gamma_t)$ is compact. But
the total oriented boundary of a compact 1-dimensional manifold is zero.

\subsection{Exercises}
\begin{itemize}

\item Try to sketch the more general proof of the orientability of the
moduli space, via the determinant line bundle of $T$: follow and adapt
the argument given in \cite{D2}.

\end{itemize}

\section{Examples and Applications}

We provide here some concrete examples in which it is possible to
compute the Seiberg--Witten invariants for some classes of four--manifolds.

\subsection{Connected Sum Theorem}

As in the case of Donaldson invariants, it is possible to prove that
Seiberg--Witten invariants vanish on connected sums. Thus the invariants
provide a test of ``irreducibility'' of the manifold.

\begin{teor}
If a manifold $X$ splits as a connected sum $X=X_1 \# X_2$, where both
pieces are manifolds with boundary with $b^{2+}(X_i)\geq 1$, then
\[ N_L (X) \equiv 0, \]
for all choices of the line bundle $L^2$.
\label{connsum}
\end{teor}

\noindent\underline{Sketch of the Proof:}There is a ``stretching the neck''
argument that works for Donaldson theory \cite{DK}, and that can be
adapted to prove the analogous result in Seiberg--Witten theory
(see the exercise in the end of the section on dimensional reduction).

\subsection{K\"ahler Manifolds}

The computation of Seiberg--Witten invariants of K\"ahler  manifolds is
described in \cite{W}; a useful review of Witten's results is also given
in \cite{FM}, where the argument is further developed for algebraic
surfaces, by algebro--geometric techniques.
See also \cite{B} for interesting
applications of Seiberg--Witten theory to K\"ahler manifolds, described in a
more algebro--geometric framework.

The first fact to mention is that a metric which is K\"ahler is non--generic,
hence the computations that we have seen in the generic case won't hold
here; in particular, the dimension of the moduli space won't be the
one given in theorem \ref{dim}.

We shall here compute the invariants under the assumption that the line
bundle $L^2$ is chosen in such a way that $c_1(L)^2=\frac{2\chi +3\sigma}
{4}$; however, we shall see that the moduli space is not just a finite
set of points.

A useful tool in the analysis of the Seiberg--Witten equations is the
following lemma:

\begin{lem}
The K\"ahler form $\omega$ induces a splitting of the spinor bundle
$S^+\otimes L$ as
\[ S^+\otimes L = (K^{1/2}\otimes L)\oplus (K^{-1/2}\otimes L), \]
where $K$ is the canonical line bundle $K=\Lambda^2(T^*X)$.
\label{K}
\end{lem}

\noindent\underline{Proof:} Consider the bundle $S$ given by
\[ S=\Lambda^{(0,*)}(X), \]
where the complexified cotangent bundle splits as $T^*(X)_{\bf C}=
T^*(X)^\prime\oplus T^*(X)^{\prime\prime}$, because of the complex structure,
and
\[ \Lambda^{(p,q)}(X)=\{ \alpha\in \Lambda^{p+q}(X)
\mid \alpha\in \Lambda^p(T^*X)^\prime\otimes
\Lambda^q(T^*X)^{\prime\prime} \}. \]
The Clifford multiplication on $S$ is defined by
\[ (\alpha^{(0,1)}+\alpha^{(1,0)})\bullet \beta=\sqrt{2}(\alpha^{(0,1)}\wedge
\beta-{\alpha^{(1,0)}}^k\beta_k), \]
where a sum over the repeated index $k$ is to be understood and the covariant
index $k$ of the 1-form $\beta$ is lowered by the metric tensor.

This spinor bundle splits as
\[ S^+=\Lambda^{(0,0)}\oplus\Lambda^{(0,2)}, \]
\[ S^-=\Lambda^{(0,1)}. \]

This will be the canonical $Spin_c$ structure associated to the
canonical line bundle $K=\Lambda^{(2,0)}$ by $S^+=1\oplus K^{-1}$.
Any other $Spin_c$ structure
is obtained by tensoring with a line bundle $L$, so that
\[ S^+=(\Lambda^{(0,0)}\oplus\Lambda^{(0,2)})\otimes L \]
and
\[ S^-=\Lambda^{(0,1)}(L). \]

\noindent QED

Note that lemma \ref{K} holds true for almost--complex structures as well:
this will be used for symplectic manifolds in the next section.

Thus we can rewrite the equations according to this splitting:
\beq
\begin{array}{c}
      {F^{(2,0)}}^+=\alpha\beta, \\
      {F^{(1,1)}}^+=-\frac{\omega}{2}(\mid\alpha\mid^2-\mid\beta\mid^2), \\
      {F^{(0,2)}}^+=\bar\alpha\bar\beta, \\
\end{array}
\label{KSW}
\eeq
where $\psi=(\alpha,-i\bar\beta)$, with $\alpha\in \Gamma(X,K^{1/2}\otimes L)$
and $\beta\in \Gamma(X,K^{-1/2}\otimes\bar L)$.

The Dirac operator has the form $D\psi=\sqrt{2}(\bar\partial\alpha
-i\bar\partial^*\bar\beta)$: full details are given in \cite{B}.

The Seiberg--Witten functional (\ref{functional}) can also be rewritten as
\beq
\begin{array}{c}
    S(A,\alpha,\beta)=\int_X (\mid F^+_A \mid^2 +<\nabla
\bar\alpha,\nabla\alpha>+<\nabla \bar\beta,\nabla\beta> \\

+\frac{1}{8}(\mid\alpha\mid^2+\mid\beta\mid^2)^2+\frac{\kappa}{4}(\mid\alpha\mid^2+\mid\beta\mid^2))dv. \\
\end{array}
\label{Kfunctional}
\eeq

Note that the functional (\ref{Kfunctional}) is symmetric under the
transformation given by
\[ \begin{array}{ccc}
         A\mapsto A, &  \alpha\mapsto\alpha, &  \beta\mapsto -\beta.
   \end{array}      \]

Hence, this has to map solutions into other solutions.
However, plugging this change of variables into the first
and third equations of (\ref{KSW}), we find that
\[ \begin{array}{cc}
         {F^{(2,0)}}^+=0, & \alpha\beta=0, \\
         {F^{(0,2)}}^+=0, & \bar\alpha\bar\beta=0. \\
   \end{array}  \]

Since $(\alpha, -i\bar\beta)$ is in the kernel of an elliptic operator
(the Dirac operator), it has some regularity properties \cite{R}, \cite{BB}.
In particular if $\alpha$ or $\beta$ vanishes on an open set, then it
vanishes identically. Hence half of the condition above reads
$\alpha\equiv 0$ or $\beta\equiv 0$.

\begin{lem}
In the Seiberg--Witten equations (\ref{KSW}), we have $\alpha\equiv 0$
or $\beta\equiv 0$, according to whether
\[ 0< \int_X \omega\wedge c_1(L) \]
or
\[ 0> \int_X \omega\wedge c_1(L). \]
\label{posneg}
\end{lem}

\noindent\underline{Proof:} Because of the decomposition of lemma
\ref{K}, $\omega$ is orthogonal to the $(0,2)$, and to the $(2,0)$
components of $F_A$. Moreover, $c_1(L)=\frac{1}{4\pi}F_A$. Hence
\[ \int_X \omega\wedge c_1(L)= -\frac{1}{8\pi}\int_X \omega\wedge\omega
(\mid\alpha\mid^2 -\mid\beta\mid^2). \]

\noindent QED

\begin{lem}
The condition ${F^{(2,0)}}^+={F^{(0,2)}}^+=0$ implies that the connection
$A$ induces a holomorphic structure on $K^{1/2}\otimes L$, or on
$K^{1/2}\otimes L^{-1}$.
\label{holom}
\end{lem}

\noindent\underline{Proof:} In holomorphic local coordinates on $X$,
a basis of self dual forms is given by the K\"ahler form
\[ \omega = \frac{i}{2}(dz_1\wedge d\bar z_1 + dz_2\wedge d\bar z_2), \]
and the forms
\[ \begin{array}{cc}
         dz_1\wedge dz_2, & d\bar z_1\wedge d\bar z_2.
   \end{array} \]
Hence ${F^{(2,0)}}^+=F^{(2,0)}$ and ${F^{(0,2)}}^+=F^{(0,2)}$. Thus
$\bar\partial_A$ gives $L$ a holomorphic structure.

\noindent QED

Suppose that $\beta\equiv 0$.
Because of lemma \ref{holom}, the equation (\ref{KSW}) becomes
\[ \bar\partial_A \alpha=0, \]
i.e. $\alpha$ is a holomorphic section of $K^{1/2}\otimes L$.

In order to compute the moduli space of solutions from these data, some
techniques of symplectic geometry are usefully applied \cite{W}.

First there is an argument which shows that we need to
consider only irreducible solutions
$\psi\neq 0$, see e.g. \cite{FM}; then the essential tools are the following
lemma and theorem.

\begin{lem}
A symplectic structure is defined on the space $\hat{\cal A}$ by:
\[ \Omega(v_A,w_A)= \int_X \omega\wedge v_A\wedge w_A \]
and
\[ \Omega(v_\alpha,w_\alpha)=\int_X \omega\wedge\omega (\bar v_\alpha
w_\alpha - \bar w_\alpha v_\alpha). \]
Here $v_A$ and $w_A$ are 1--forms in the tangent space to the connection $A$,
and $v_\alpha$ and $w_\alpha$ are in the tangent space to the section $\alpha$.
\end{lem}

\begin{teor}
The gauge group ${\cal G}$ acts on this symplectic manifold. By a
reduction argument it is possible to identify the moduli space $M$ with
the quotient of the fibre over zero of the moment map by the action of
the complexified gauge group ${\cal G}_{\bf C}={\cal M}(X,{\bf C}^*)$:
\[ M= \mu^{-1}(0)/{\cal G}_{\bf C}. \]
\label{Kmod}
\end{teor}

\begin{corol}
As a consequence of theorem \ref{Kmod}, the moduli space of solutions of the
Seiberg--Witten equations on a K\"ahler manifold is given by the set of
possible choices of a holomorphic structure on $K^{1/2}\otimes L$ and
the corresponding projectivisation of the space of holomorphic sections,
\[ {\bf P}H^0 (X, K^{1/2}\otimes L). \]
\end{corol}

\noindent\underline{Proof of Corollary:} The same holomorphic structure on
$K^{1/2}\otimes L$ is determined by $A$ and by $\lambda^{-1}
A\lambda$, for $\lambda \in {\cal G}_{\bf C}$. A section $\alpha$ of
$K^{1/2}\otimes L$ is acted upon by the complexified gauge group via
$\alpha\mapsto\lambda\alpha$. Hence modding out by the action we get the
projectivisation.

\noindent QED

If $\alpha=0$ and $\beta\neq 0$, just switch $K^{1/2}\otimes L$ with
$K^{1/2}\otimes L^{-1}$ in the above result.

\subsection{Symplectic Manifolds}

The same technique used above in attacking the problem of K\"ahler
manifolds can be partially extended to symplectic manifolds.

The following theorem is proven in \cite{T}:
\begin{teor}
If $X$ has a symplectic structure $\omega$ compatible with the orientation,
and it has $b^{2+}>1$, then the
Seiberg--Witten invariant corresponding to the canonical $Spin_c$--structure
(see lemma \ref{K}) is
\[ N_K(X)=\pm 1. \]
\end{teor}

The proof of this theorem given in \cite{T}, has subsequently been
considerably simplified \cite{T6}, as reported also in \cite{T2}.

A brief sketch of the argument is as follows.

\noindent\underline{Sketch of the proof:} The main idea is to show that
the Seiberg--Witten invariant can be computed by means of a one parameter
family of perturbed equations.

These are of the form
\[ D_A\psi =0 \]
and
\[ F^+_A=F^+_{A_0}+\frac{1}{4}<e_ie_j\psi,\psi>e^i\wedge e^j -\frac{ir}{4}
\omega. \]

One particular solution $(A_0,u_0)$, that corresponds to the value $r=0$
of the parameter, is constructed by taking
$A_0$ to be a connection on $K^{-1}$, and
the section $u_0$ to be a covariantly constant norm 1 section of the
trivial summand in $S^+=1\oplus K^{-1}$ (lemma \ref{K}).
In the condition $\nabla_{A_0}u_0=0$ we are considering
the induced covariant derivative on
$S^+=1\oplus K^{-1}$, obtained from the exterior derivative $d$, and
$\nabla_{A_0}=\frac{1}{2}(1+i\omega)
\tilde\nabla_{A_0}$, with the notation of definition \ref{twist}.

In the perturbed equation above the section $\psi$ is decomposed as
\[ \psi=r^{1/2}(\alpha u_0 +\beta), \]
with $\alpha$ a function, and $\beta$ a section of $K^{-1}$.

By successively projecting the second equation on the trivial summand and on
$K^{-1}$, and by Sobolev techniques, a uniform bound (independent of $r$) is
obtained on the expression
\[ \int_X (\frac{4}{r}\mid F^+_a \mid^2 +\mid\nabla_a\mid^2 -\frac{r}{8}
(1-\mid\alpha\mid^2)^2)dv, \]
where $F_a=\frac{1}{2}(F_A-F_{A_0})$.

Hence considering a sequence of solutions $(A_m,\psi_m)$, corresponding to
a sequence of parameters $r_m\to\infty$, the above says that
$\mid\alpha_m\mid\to 1$, and $\mid\nabla_a\alpha_m\mid\to 0$.

Finally it is shown that the $\beta_m$ vanish identically for large
enough $m$; and consequently $\alpha_m$ is covariantly constant of
norm 1.

The upshot is that $(A_m,\psi_m)$, for large enough $m$, is gauge equivalent
to $(A_0,u_0)$; and therefore there is just one point in the
moduli space, and the Seiberg--Witten invariant is $\pm 1$.

\noindent QED

As a consequence of Taubes' result that symplectic manifolds admit some
non--trivial Seiberg--Witten invariants, it is natural to ask whether a
converse might hold as well, i.e. whether any four-manifold with non-trivial
invariants might be symplectic. This problem has been considered in
\cite{KMT}, and solved negatively. There are 4-manifolds with non--trivial
Seiberg--Witten invariants but without symplectic structure.

The strategy used to construct such examples is to prove two lemmata:
\begin{lem}
If $X$ is symplectic and it decomposes as a smooth connected sum
$X=X_1\# X_2$, then one of the summands, say $X_1$, has negative
definite intersection form (because of the connected sum theorem)
and fundamental group $\pi_1(X_1)$ which doesn't admit any non-trivial
finite quotient.
\label{l1}
\end{lem}

\begin{lem}
If $X_1$ has some non-trivial Seiberg--Witten invariant and $X_2$ has
$b^1(X_2)=b^{2+}(X_2)=0$, then the connected sum $X=X_1\# X_2$ also
has some non-trivial Seiberg--Witten invariant.
\label{l2}
\end{lem}

The first lemma depends on the connected sum theorem, and on the behaviour
of $b^{2+}$ on finite covers; the second depends on an explicit construction
of a $Spin_c$ structure on $X$, given the structure on $X_1$ with
non-trivial invariants and a unique $Spin_c$ structure on $X_2$.
A ``stretching the neck'' argument is required to glue solutions
corresponding to the two separate structures, when performing the connected
sum, \cite{Wa}.

Lemma \ref{l1} is further refined in \cite{K}, where it is also proven that,
under the same hypotheses, $X_1$ is an integral homology sphere.

Manifolds with non-trivial invariants and with no
symplectic structures are therefore given by the following class of
examples.

\begin{esem}
Let $X_1$ be a symplectic manifold with $b^{2+}(X_1)>1$; by \cite{T}
$X_1$ has a non-trivial invariant. Let $X_2$ be a manifold with
$b^1(X_2)=b^{2+}(X_2)=0$, and with fundamental group that admits
non-trivial finite quotients. Then $X=X_1\# X_2$ has a non-trivial invariant,
but it doesn't admit a symplectic structure.
\end{esem}

In \cite{K} another interesting problem is considered. After the results of
\cite{T}, \cite{T1}, it was conjectured by Taubes that the following
decomposition result should be true:

\begin{conj}
Every smooth, compact, oriented 4--manifold without boundary is a
connected sum of symplectic manifolds, with either the symplectic or
the opposite orientation, and of manifolds with definite intersection form.
\end{conj}

In \cite{K} it is proven that this conjecture is false. The construction
is similar to the arguments used in the two lemmata \ref{l1},\ref{l2}.

Another conjecture remains open at the moment:
\begin{conj}
Every smooth, compact, oriented, simply connected four--manifold
without
boundary is a connected sum of symplectic manifolds, with both
orientations allowed.
\end{conj}

Note that the latter conjecture would imply the smooth Poincar\'e conjecture.

\subsection{The Thom Conjecture}

One of the first applications of the Seiberg--Witten gauge theory was
to prove the Thom conjecture, for ${\bf C}P^2$ in \cite{KM1}; a more general
version of this result will be presented in an upcoming paper by Morgan,
Szabo, and Taubes \cite{MST}.

We'll reproduce here the proof that is given in \cite{KM1}.

\begin{teor}
An oriented two--manifold $\Sigma$ that is embedded in ${\bf C}P^2$ and
represents the same homology class as an algebraic curve of degree $d>3$,
has genus $g$ such that
\[ g \geq \frac{(d-1)(d-2)}{2}. \]
\label{thom}
\end{teor}

\noindent\underline{Sketch of the proof:} The proof is obtained in several
steps.

(1) good metrics:
\begin{defin}
a metric on a manifold $X$ is ``good'' with respect to a certian choice
of the line bundle $L^2$ if there is no non-trivial singular
point in the moduli space. Lemma \ref{poscurv} implies that this happens
whenever $c_1(L)$ is not in $H^{2-}(X;{\bf R})$.
\label{good}
\end{defin}

\noindent We shall consider the manifold
\[ X={\bf C}P^2 \# n\overline{{\bf C}P^2}, \]
the blow-up (in the language of algebraic geometry) of the projective space
at $n$ points, and identify the good metrics on $X$ in terms of a condition
on the product of the first Chern class of the line bundle
with a harmonic form that depends on the metric:
\[ \int_X c_1(L)\cup [\omega_g]\neq 0. \]

(2) If the product $\int_X c_1(L)\cup [\omega_g]$ is negative, the
Seiberg--Witten moduli space is empty.

This is shown by first proving that the Seiberg--Witten invariant (mod 2)
changes
parity if $\int_X c_1(L)\cup [\omega_g]$ changes sign; and then proving
that there is a particular choice of the metric $g$ such that
$\int_X c_1(L)\cup [\omega_g]>0$ and the moduli space is empty.

(3) Suppose a four--manifold $X$ splits along a three manifold $Y$, so that
the metric is a product on a neighbourhood $[-\epsilon,\epsilon ]\times Y$.
Consider the metric $g_R$ given by inserting a flat cylinder $[-R,R]\times Y$.
If the moduli space $M(g_R)$ is non-empty for all large $R$, then there
exists a solution of the Seiberg--Witten equations which is ``translation
invariant
in a temporal gauge'' on the manifold ${\bf R}\times Y$. This means that
the $dt$ component of the connection $A$ vanishes (see definition
\ref{tempgau} below).

(4) If there is a solution on ${\bf R}\times Y$ that is translation
invariant in a temporal gauge and $Y=S^1\times \Sigma$, where $\Sigma$ is a
surface
of constant scalar curvature and genus $g\geq 1$, then there is an estimate
\[ \mid \int_\Sigma c_1(L^2) \mid < 2g-2. \]

(5) Let $H$ be the generator of $H^2({\bf C}P^2;{\bf Z})$; by assumption,
we have an embedded
surface $\Sigma$ of genus $g$ that determines the homology class
dual to $dH$.

Consider the embedding
\[ \Sigma \hookrightarrow {\bf C}P^2 \# d^2 \overline{{\bf C}P^2}. \]
The second cohomology of $X={\bf C}P^2 \# d^2 \overline{{\bf C}P^2}$ has
generators $H$ and $E_i$, $i=1,\ldots d^2$, with intersection form
$Q_X=(1,d^2)$.

Take $\tilde\Sigma=\Sigma \# d^2 S^2$, where  $S^2 \in \overline{{\bf C}P^2}$
is dual to $-E_i$. Thus the homology class $[\tilde \Sigma ]$ is dual to
$dH - E$, $E=\sum_i E_i$.

Take a tubular neighbourhood $T$ of $\tilde\Sigma$ and a metric $g_0$ on $X$
such that $Y=\partial T =\tilde\Sigma \times S^1$ with a product metric and
constant scalar curvature $-2\pi (4g-4)$ on $\tilde\Sigma$ (assume $\tilde
\Sigma$ has unit area and use Gauss--Bonnet).

Consider the canonical line bundle $L^2=K$. This has Chern class
$c_1(K)=3H-E$.

Insert a cylinder $[-R,R]\times Y$. Then
\[ \int_X c_1(K)\cup [\omega_{g_R}]= [\tilde\Sigma ] [\omega_{g_R}] -
(d-3)\int_X H\cup [\omega_{g_R}]. \]

By normalising $1=\int_X H\cup [\omega_{g_R}]$, and showing that
$[\tilde \Sigma ][\omega_{g_R}] \to 0$ as $R\to\infty$, the result of
the theorem follows from the estimate of step (4).

Now we'll see in more details the various steps of the proof. Part of it
is left as a series of exercises at the end of this section.

\noindent\underline{Step (1):} The intersection form $Q_X$ on the manifold
$X={\bf C}P^2 \# n\overline{{\bf C}P^2}$ has signature $(1,n)$; thus it
determines a cone $C$ in $H^2(X;{\bf R})$ where the form is positive.

If $H$ is the generator of the cohomology ring of ${\bf C}P^2$, then
$H\in C$. Call $C^+$ the fold of the cone that contains $H$.

The manifold $X$ has $b^{2+}=1$; hence for a chosen metric there exists a
unique harmonic self dual form $\omega_g$ such that the corresponding
cohomology class $[\omega_g ]\in C^+$.

So we have that the metric $g$ is good in the sense of lemma \ref{good} iff
\[ \int_X c_1(L)\wedge \omega_g \neq 0. \]

\noindent\underline{Step (2):} $\int_X c_1(L)\wedge \omega_g = 0$ detects the
presence of some singular point in the moduli space. By perturbing the
equation with a small $\eta\in \Lambda^{2+}$, we can assume that on a
given path of metrics $\{ g_t \mid t\in [0,1] \}$ the expression
\[ f(t)=\int_X C_1(L)\wedge \omega_{g_t} +
2\pi \int_X \eta\wedge \omega_{g_t} \]
changes sign transversally at $t=0$.

Hence the parametrized moduli spaces $M_{g_t}$ look like
a family of arcs; we want to prove that at a singular point an odd number of
arcs meet. Thus the invariant computed mod 2 changes parity.

This requires the analysis of a local model of the moduli space around a
singular point. A model of Donaldson's can be adapted to this case
\cite{D3}.

\noindent\underline{Step (3):} This is the part of the proof where the
gauge theoretic techniques have a prominent role. Since the argument
involves the dimensional reduction of Seiberg--Witten theory to three
dimensions, we postpone the proof of this step until after the section that
deals with three--manifolds applications of the theory.

\noindent\underline{Step (4):} Assuming $\Sigma$ to be of unit area and with
constant scalar curvature, the Gauss--Bonnet theorem implies that the scalar
curvature $\kappa = -4\pi (2g-2)$. From the estimate on the spinor $\psi$
given in lemma \ref{psibound}, we have
\[ \mid \psi \mid^2 \leq 4\pi (2g-2). \]

But, from the equation (\ref{eqSW2}), since
\[ \mid <e_i e_j\psi,\psi>\mid^2 =2\mid \psi \mid^4, \]
we get an estimate $\mid F^+_A \mid \leq \sqrt{2} \pi (2g-2)$.

Since the solution is translation invariant in a temporal gauge, $F_A$
is the pullback on ${\bf R}\times Y$ of a form on $Y$. Hence $F_A\wedge
F_A =0$, and this means that $\mid F_A^+ \mid =\mid F_A^- \mid$, since
$F_A\wedge F_A=(\mid F_A^+ \mid^2 - \mid F_A^- \mid^2)dv$.

Thus the resulting estimate on $F_A$ is
\beq
\label{est}
\mid F_A \mid \leq 2\pi (2g-2).
\eeq
Therefore
\[ \mid \frac{1}{2\pi}\int_\Sigma F_A \mid = \mid c_1(L^2)[\Sigma] \mid
\leq 2g-2. \]

\subsection{Other Applications}

Here we want to present briefly some of the current directions of research
that give interesting applications of the Seiberg--Witten theory.

In the field of differential geometry it is worth mentioning some results on
Einstein metrics. A uniqueness theorem for such metrics on compact quotients
of irreducible 4-dimensional symmetric spaces of non-compact type is proven
in \cite{LeB}.

In the study of K\"ahler metrics, besides the papers mentioned
in a previous section, we ought to mention \cite{LeB2} on extremal K\"ahler
metrics, \cite{Leu}, and \cite{OT}.

Another approach to Seiberg--Witten theory is related to integrable
systems; so far it has been mainly developed in Physics papers like
\cite{GKMMM}.

It is certainly worth mentioning the contribution of \cite{B} where,
in the case of K\"ahler manifolds, the language of Fulton's intersection
theory is adapted to a Fredholm context in such a way as to extend the
definition
of the invariants, with no need of the generic condition or of perturbation.

\subsection{Exercises}
\begin{itemize}

\item Fill in the details of the proof of step (1) above.

\item Check that the argument given in \cite{D3} can be adapted to the
proof of step (2).

\item To complete step (5): check the details of the proof that
$[\tilde \Sigma ][\omega_{g_R}] \to 0$ as $R\to\infty$, as given in
\cite{KM1}.

\end{itemize}

\section{Quantum Field Theory}

What we would like to introduce here is an application of the Mathai--Quillen
formalism \cite{MQ} to the Seiberg--Witten gauge theory.
The purpose is to construct a Lagrangian,
given the moduli space.

\subsection{A Topological Lagrangian}

The references available for this section are \cite{LabM}, \cite{ZWC},
\cite{GMR}. The main reference about Quantum Field Theory and
four-manifolds invariants is Witten's paper \cite{W2}.

We'll try here to enlighten mainly the
mathematical aspects. To this purpose, we shall follow the very clear
and enlightening introduction to the Mathai--Quillen formalism and its
applications to Quantum Field Theory, given in \cite{AJ}.

The purpose of the QFT approach is to construct a topological Lagrangian,
in the sense of Witten \cite{W2}, for our gauge theory such
that the Seiberg--Witten invariants can be obtained as the expectation values
of certain operators. In a more mathematical language, this
corresponds to an infinite dimensional analogue (\cite{MQ}, \cite{AJ}) of
the Chern--Gauss--Bonnet theorem.

In the finite dimensional case, the Euler number of a $2m$--dimensional
vector bundle $E$ over a manifold $Y$ of dimension $2m$ is obtained
by evaluating on $Y$ a certain class $\omega$: $e=\int_Y \omega$.

A different way to write the integral above is given in \cite{MQ}: as
$\int_Y \sigma^*U$, where $U$ is a form on $E$ representing the Thom class,
with Gaussian decay along the fibres $E_x$;
$\sigma$ is a section of $E$.

In case $n=dim(Y)>2m$, the integral above is replaced by $\int_Y
\eta\wedge\omega$,
where $\eta$ is an element of the $(n-2m)$--cohomology of $Y$.

In \cite{MQ} an explicit form for a representative of the Thom class is
constructed as follows.

\begin{defin}
Given a vector bundle $E=P\times_G V$ obtained from a principal $G$
bundle $P$ via a representation
$\rho:G\to SO(2m)$, the Mathai--Quillen form is
\[ U=\pi^{-m} e^{-x^2} \int \exp (\frac{w^t\Omega w}{4}+i dx^t\cdot w)
{\cal D}w, \]
where the $w$'s are a basis of the tangent bundle to the fibre, $\Omega$
is the image of the curvature matrix of $G$ under the representation $\rho$,
and the symbol ${\cal D}w$ means that, in the power series expansion,
only the coefficient of $w_1\wedge\cdots\wedge w_{2m}$ is to be taken.
\label{MaQu}
\end{defin}

Some manipulations of the above expression, that are fully explained in
\cite{AJ}, lead to a formula for the Euler number of the bundle $E$.

\begin{teor}
The Euler class of the bundle ${\cal E}$ is computed in terms of the
Mathai--Quillen form as the integral over $P$ of the form
\beq
\begin{array}{c}
2^{-d}\pi^{-d-m}\int \exp (-\mid \sigma \mid^2+\frac{w^t\rho(q)w}
{4} +i d\sigma^t w \\
-i<d\nu,h>+i(q,Rh)+<dy,Cf>)
{\cal D}f {\cal D}w {\cal D}q {\cal D}h.
\end{array}
\label{MQform}
\eeq
\label{Eunum}
\end{teor}

In the above expression $d=dim(G)$; $q$, $h$, and $f$ are
Lie algebra variables that arise in the use of a Fourier transform and
in expressing the invariant volume in terms of the Killing form;
the $y$'s are coordinates on $P$; and $\nu$ is a canonical 1--form
with values in ${\cal L}(G)$, defined by the action of $G$ on $P$.
$C$ and $R$ arise in the definition of $\nu$ as
\[ \nu_\xi (v)=<C\xi,v>, \]
where $v$ is a tangent vector, and $\xi\in {\cal L}(G)$; $C\xi$ is the
vector field given by the infinitesimal action of $\xi$ (like in definition
\ref{hamilt}). The operator $R$ is defined as $R=C^*C$, i.e.
$<C\xi,C\zeta>=(R\xi,\zeta)$, where the inner product $(,)$ is the one
given by the Killing form.

The argument given here carries over, at least formally, to some infinite
dimensional cases, \cite{MQ}, \cite{AJ}, and it therefore provides the right
mathematical setup in
which the topological Lagrangian introduced by Witten for Donaldson
theory \cite{W2} lives. The analogous construction works for
Seiberg--Witten theory, as shown in \cite{LabM}, or \cite{ZWC}.

Consider the space $\hat{\cal A}$ of $U(1)$--connections and non-zero spinors.
We have already seen in (\ref{chcomplex}) that the tangent space of
$\hat{\cal A}$ at a point $(A,\psi)$ is given by
\[ \Lambda^1(X)\oplus \Gamma (X,S^+\otimes L). \]

Consider a trivial bundle ${\cal E}$ over $\hat{\cal A}$ with fibre
\[ {\cal F}=\Lambda^{2+}\oplus \Gamma(X,S^- \otimes L). \]

\begin{lem}
The set of irreducible solutions of the Seiberg--Witten
equations can be described as the zero set of a section of ${\cal E}$
given by
\[ \sigma(A,\psi)=(F^+_A-\frac{1}{4}<e_ie_j\psi,\psi >e^i\wedge e^j,
D_A\psi). \]
\label{quant1}
\end{lem}

This is clear from the construction. Moreover, since the gauge group acts
freely on this space of solutions by (\ref{action}), we can induce a section
\[ \bar\sigma: \hat{\cal A}/{\cal G}\rightarrow {\cal E}/{\cal G}, \]
which describes the smooth part of the moduli space.

\begin{lem}
The section $\sigma(A,\psi)$ satisfies
\[ \mid \sigma(A,\psi) \mid^2= S(A,\psi), \]
where $S(A,\psi)$ is the Seiberg--Witten functional defined in
(\ref{functional}).
\label{qsect}
\end{lem}

\noindent\underline{Proof:}This follows by direct computation from the
definition of the $L^2$--inner product of forms and sections in
\[ {\cal F}=\Lambda^{2+}\oplus \Gamma(X,S^- \otimes L). \]
Further details are left as an exercise.

\noindent QED

Now we can identify the various terms of the Euler class, as given in
theorem \ref{Eunum}: this has been done in \cite{LabM}, \cite{ZWC}.

Note that, from (\ref{chcomplex}), we shall have variables $dy=(\alpha,\phi)$,
with $\alpha\in \Lambda^1(X)$ and $\phi\in \Gamma(X,S^+\otimes L)$,
that represent a basis of forms on $\hat{\cal A}$; and the variable $f$
which counts the ``gauge directions'': $f\in \Lambda^0(X)$ satisfies
$e^{if}\in {\cal G}$. Similarly $q$ and $h$ are in the Lie algebra, i.e.
in $\Lambda^0(X)$. $w=(\beta,\chi)\in \Lambda^{2+}(X)\oplus
\Gamma (X,S^-\otimes L)$ is the variable along the fibre.

The term $d\bar\sigma$, therefore, is the linearization of the
Seiberg--Witten equations given in lemma \ref{linear}:
\[ d\bar\sigma(A,\psi)=(D_{A_0}\phi + i\alpha\cdot \psi_0,
 d^+\alpha -\frac{1}{2}Im(<e_ie_j\psi_0,\phi>)e^i\wedge e^j). \]

Moreover, it is clear that the operator $C$ is the map $G$ of the complex
(\ref{chcomplex}) that describes the infinitesimal action of the gauge
group. Hence
\[ C(f)=G(f)=(-idf, if\psi)\in \Lambda^1(X)\oplus\Gamma(X,S^+\otimes L) \]
and
\[ C^*(\alpha,\phi)=-d^*\alpha +\frac{1}{2}Im(<e_ie_j\psi,\phi>). \]
Thus the operator $R$ is given by
\[ R=d^*d+\mid\psi\mid^2. \]

Thus all the terms in (\ref{MQform}) can be computed:
\[ \frac{1}{4}w^t\rho(q)w =\frac{-i}{4} q
\mid\chi\mid^2 dv, \]
\[ i d\bar\sigma^t w =i\beta\wedge *d^+\alpha
-\frac{i}{2}(\beta,Im(<\psi,\phi>)e^i\wedge e^j)dv +
<D_A\phi +i\alpha\psi,\chi>dv, \]
\[ -i<d\nu,h>=-ih\mid\beta\mid^2, \]
\[ i(q,Rh)=iq(d^*dh+\mid\psi\mid^2 h), \]
\[ <dy,Cf>=(C^*dy,f)=df\wedge *d\alpha +\frac{f}{2}Im(<e_ie_j\psi,\phi>)dv.\]

The only expression that requires some more comments is the third above,
where $dC^*(dy_1,dy_2)$ is computed \cite{LabM} using
\[ dC^*(dy_1,dy_2)=dy_1(C^*(dy_2))-dy_2(C^*(dy_1))-C^*([dy_1,dy_2]) \]
and constant vector fields $dy_1$ and $dy_2$, whose Lie bracket vanishes
identically.

Thus we have constructed a ``topological Lagrangian''; we want to recover
the Seiberg--Witten invariants as correlation functions.

\subsection{Seiberg--Witten Invariants Revisited}

There is a ``natural'' definition of the Seiberg--Witten
invariants when the moduli space $M_L$ is zero-dimensional, which
we introduced in definition \ref{inv0}. As already mentioned, however,
the definition of the invariants for a higher dimensional moduli space
implies a certain arbitrariness in the choice of the cohomology class to
pair with the cycle given by the moduli space itself.

We shall see here how to get a definition like \ref{invpos} by a different
approach. This should reinforce the feeling that there ought to be a
``best'' choice that generalises definition \ref{inv0}.

What follows is partially just a formal argument, as the computations have
to be carried out in an infinite dimensional setup.

In the previous paragraph we have shown that
the Mathai--Quillen form of the Euler class of the bundle $E$ on
$\hat{\cal A}/{\cal G}$ is
\[ e= 2^{-d}\pi^{-m-d}\int \exp (-\mid \sigma \mid^2+(\frac{-i}{4}q\mid\chi
\mid^2-\frac{i}{2}(\beta,Im<\psi,\phi>e^i\wedge e^j) \]
\[ +<D_A\phi+i\alpha\psi,
\chi>-ih\mid\beta\mid^2+iq(\Delta h+\mid\psi\mid^2h)+\frac{f}{2}Im<e_i
e_j\psi,\phi>)dv \]
\[ +i\beta\wedge*d^+\alpha+df\wedge*d\alpha) {\cal D}f {\cal D}h
{\cal D}q {\cal D}\beta {\cal D}\chi. \]

When the moduli space is zero dimensional the Seiberg--Witten invariant
is obtained as the Euler number of the bundle, \cite{AJ},
i.e. by integrating the
above class over the total space $P$ (here ${\cal D}\alpha {\cal D}\phi$
is just a formal measure, or a ``functional integration'' in the language
of physics, since the form is being integrated over an infinite dimensional
manifold).

The Euler class above formally has ``codimension'' equal to the dimension
of the moduli space, since that is $-Ind(C^*)$ where $C^*$ is the
chain complex obtained by linearizing the Seiberg--Witten equations
(\cite{AJ}).

Hence, when the moduli space is of positive dimension, in order to obtain
numerical invariants, we need to cup the class above with other
cohomology classes of $\hat{\cal A}/{\cal G}$.

The following is a possible construction, mimicking Donaldson's construction
of the polynomial invariants. The idea is to choose a bundle over
$\hat{\cal A}/{\cal G}\times X$, integrate a characteristic class of
this bundle against a homology class of $X$, and then, restricting to
$M\hookrightarrow \hat{\cal A}/{\cal G}\times X$, get a map
\[ \omega : H_k(X)\rightarrow H^{i-k}(\hat{\cal A}/{\cal G}), \]
where $i-k$ is the dimension of $M$, and $i$ is the degree of the
characteristic class.

There is a $U(1)$-bundle $\hat{\cal A}/{\cal G}_0\rightarrow
\hat{\cal A}/{\cal G}$, where ${\cal G}_0$ is the gauge group of base point
preserving maps. In fact ${\cal G}/{\cal G}_0=U(1)$.  The restriction of
this bundle to $M\hookrightarrow \hat{\cal A}/{\cal G}$ is the $U(1)$
bundle used in definition \ref{invpos}.

We can form the bundle ${\cal Q}=\hat{\cal A}\times_{{\cal G}_0} L^2$ over
$\hat{\cal A}/{\cal G}\times X$ and consider its first Chern class
$c_1({\cal Q})$.

Then take the pullback via the inclusion of the moduli space
$M\hookrightarrow \hat{\cal A}/{\cal G}$; this gives $i^*{\cal Q}\rightarrow
M\times X$.

The construction above gives a map
\[ \omega : H_k(X;{\bf Z})\rightarrow H^{2-k}(M;{\bf Z}). \]
In fact, take a class $\alpha\in H_k(X;{\bf Z})$, and its Poincar\'e dual
$a=PD(\alpha)$; take the first Chern class
\[ c_1({\cal Q})\in H^2(\hat{\cal A}/{\cal G}\times X;{\bf Z}) \]
and decompose it according to the K\"unneth formula
\[ H^2(\hat{\cal A}/{\cal G}\times X;{\bf Z})=\oplus_j
H^j(X;{\bf Z})\otimes H^{2-j}(\hat{\cal A}/{\cal G};{\bf Z}), \]
as
\[ c_1({\cal Q})=\oplus_j c_1(Q)^j_{2-j}. \]

Evaluate against $X$ to get a class
\[ \int_X c_1({\cal Q})^k_{2-k}\wedge a \in
H^{2-k}(\hat{\cal A}/{\cal G};{\bf Z}). \]
The pullback via $i^*$ defines a class in $H^2(M;{\bf Z})$:
\[ \omega(\alpha)\equiv i^* \int_X c_1({\cal Q})\wedge a. \]

This defines maps
\[ q_d:H_{k_1}(X)\times\cdots\times H_{k_r}(X)\rightarrow {\bf Z} \]
with $\sum_{j=1}^{r} (2-k_j) =d$
\[ q_d(\alpha_1,\cdots,\alpha_r)=\int_{M} \omega(\alpha_1)\wedge
\cdots\wedge\omega(\alpha_r), \]
where $d$ is the dimension of the moduli space $M$. These play the role,
in our construction, of Donaldson's polynomial invariants.

Thus, according to the quantum field theoretic formalism, we obtain
the invariants by evaluating over $\hat{\cal A}/{\cal G}$ the Euler
class cupped with classes $\omega(\alpha)$:
\[ N_L\equiv \int_{\hat{\cal A}/{\cal G}} \bar\sigma^*(e)\wedge\omega(\alpha_1)
\wedge\cdots\wedge\omega(\alpha_r), \]
with $d=dim(M_L)$.

This should be the definition that corresponds to the
description of the invariants as
expectation values of the operators obtained by the formalism
of Quantum Field Theory, \cite{W2}.

The operators constructed in \cite{LabM}, or \cite{ZWC} following \cite{W2},
are:
\[ W_{k,0}=\frac{f^k}{k!}, \]
\[ W_{k,1}=\alpha W_{k-1,0}, \]
\[ W_{k,2}=F W_{k-1,0}-\frac{1}{2} \alpha\wedge\alpha W_{k-2,0}, \]
\[ W_{k,3}=F\wedge\alpha W_{k-2,0}-\frac{1}{3!} \alpha\wedge\alpha\wedge\alpha
W_{k-3,0}, \]
\[ W_{k,4}=\frac{1}{2} F\wedge F W_{k-2,0} -\frac{1}{2}
F\wedge\alpha\wedge\alpha W_{k-3,0}
-\frac{1}{4!} \alpha\wedge\alpha\wedge\alpha\wedge\alpha W_{k-4,0}. \]

The choice of different $k$ should correspond \cite{W2} to different
choices of characteristic classes to pair with the homology
classes of $X$ in the construction of the invariants.
In particular we should obtain
a relation between some of these operators and the polynomial invariants
constructed above,
adapting  to the present case the argument given for Donaldson theory
in  \cite{AJ}.

The Chern class $c_1({\cal Q})$ should be interpreted as
a curvature on the infinite dimensional bundle ${\cal Q}$ in such a way that
the components $c_1({\cal Q})^{2-i}_i\in \Lambda^i(X)\otimes \Lambda^{2-i}
(\hat{\cal A}/{\cal G})$ should be written in terms of the operators $W$ as
\[ \int_{\hat{\cal A}/{\cal G}}\bar\sigma^*(e)
\prod_{i=1}^r\int_{\alpha_i} W_i =\int_M \int_{\alpha_{1}}
c_1({\cal Q})^{2-i_1}_{i_1}\wedge\cdots\wedge\int_{\alpha_r}
c_1({\cal Q})^{2-i_r}_{i_r}, \]
integrated over submanifolds $\alpha_i$ of $X$ of the proper dimension.

It is not clear that the polynomial invariants defined in this section are
non-trivial for positive dimensional moduli spaces.
In fact it is conjectured (\cite{Au} pg.34) that for simple type
manifolds (see the section on Seiberg-Witten and Donaldson theory) the
only non-trivial invariants are associated to zero-dimensional moduli
spaces.

\subsection{Another Construction}

The main problem with the above argument is to somehow make precise the
construction of homology and cohomology classes of infinite dimensional
Hilbert manifolds.

It is worth mentioning that
the same problem is dealt with in \cite{B}, with a different approach.
R. Brussee adapts to the Seiberg--Witten case a construction of
Pidstrigatch \cite{Pi}, \cite{PiT} that defines a fundamental cycle
in $H_d(M)$, where $d$ is the generic dimension of a moduli space
$M$ arising from an elliptic equation.
In order to carry out his construction, Brussee considers the case
of K\"ahler manifolds, where he can use the formalism of algebraic geometry
(intersection theory).

The cycle represents the
homological Euler class of an infinite dimensional bundle over $M$.
Hence, we can consider this approach as a way to make precise the idea that
the Mathai-Quillen-Euler class ought to have codimension $d$.

\subsection{Exercises}
\begin{itemize}

\item Complete the computation needed for the proof of lemma \ref{qsect}.

\item This and the following exercises are intended to address some interesting
questions that arise from the QFT formalism introduced above. The first
problem is to describe more precisely the relation between the
infinite dimensional bundle
used to compute the invariants in this context, and the choice made
by \cite{T2} (definition \ref{invpos}).
The two following problems deal with supersymmetry.

\item The formalism above that produces the topological Lagrangian can
be rewritten in terms of superalgebras, \cite{W2}.
Try to follow the argument given in \cite{LabM}, \cite{ZWC}.

\item We already know that the Seiberg--Witten functional provides many other
critical points that are non--minimising. Therefore, they do not
correspond to solutions of the Seiberg--Witten equations, but rather of the
second order variational problem (\ref{var1}), (\ref{var2}).
Thus, we would like to modify the functional (\ref{functional}) in such a way
that it still encodes all the information concerning the solutions of the
Seiberg--Witten equations, but also in such a way to get rid of all
non--minimising critical points.
In Physics this kind of problem is taken care of exactly in the supersymmetric
formulation of gauge theories. How does this relate to the results of this
section?

\end{itemize}

\section{Dimensional Reduction and Three Manifolds}

So far we have always considered Seiberg--Witten gauge theory as a tool
for studying the differentiable structure of four-manifolds. In the present
section we try a dimensional reduction, that produces a three-dimensional
gauge theory.

At the moment, not many three-manifold results are known, but it seems
that the analogues of the four-dimensional invariants provide an
invariant of three manifolds which resembles in nature the Casson
invariant. This at least is known at the level of Physics, in a quantum
field theoretic formalism, \cite{ZWC}; but the problem is
still to be investigated on a rigorous mathematical ground.
Certainly this will be done in a very short time, \cite{T3}; anyhow, here
we'll try to present what is known at the moment, and what is expectable.
As a reference for this section the reader can look at part of
\cite{KM1}, and at \cite{ZWC}.

Consider a manifold $X=[0,1]\times Y$, where $Y$ is a compact connected
three-manifold without boundary, endowed with a cylindrical metric.

$X$ is a four-manifold with boundary. Although we have so far been considering
the Seiberg--Witten equations on four-manifolds without boundary,
we may as well define the same equations on $X$. We have to choose a
$Spin_c$-structure; and we can do so in such a way that the positive and
negative spinor bundles $S^\pm \otimes L$ are the pullback, via the
projection on the second factor, of a $U(2)$ bundle over $Y$, $\tilde S
\otimes L$, which endows $Y$ with a $Spin_c$-structure.

The positive and negative spinor bundles over $X$ are isomorphic, via Clifford
multiplication by $dt$, $t\in [0,1]$.

\begin{defin}
We say that a pair $(A,\psi)$ on $X$ is {\em in a temporal gauge} if
the $dt$ component of the connection $A$ is identically zero.
\label{tempgau}
\end{defin}

A pair $(A,\psi)$ on $X$ that is in a temporal gauge induces a path
of connections and sections of the spinor bundle over $Y$ by defining
$A(t)$ and $\psi(t)$ to be the restrictions of $A$ and $\psi$ to
$\{ t \}\times Y$.

We can perform a dimensional reduction as follows.
\begin{teor}
For a pair $(A,\psi)$ in a temporal gauge, the Seiberg--Witten equations
(\ref{eqSW1}), (\ref{eqSW2}) on $X$
induce the following equations on $Y$:
\beq
\frac{d}{dt}\psi=-\partial_A\psi
\label{3eqSW1}
\eeq
and
\beq
\frac{d}{dt}A=-*F_A+<e_i\psi,\psi>e^i.
\label{3eqSW2}
\eeq
\label{dimred}
\end{teor}

\noindent\underline{Sketch of the proof:} The Dirac operator on $X$ twisted
with the connection $A$ has the form
\[ D_A=\left( \begin{array}{cc}
                 0 & D^+_A \\
                 D^-_A & 0
              \end{array} \right), \]
\[ D^+_A = \partial_t +\partial_{A(t)}, \]
where $\partial_A$ is the self-adjoint Dirac operator on $Y$ twisted with
the time dependent connection $A(t)$.

For the curvature equation (\ref{3eqSW2}), observe that 1-forms on $Y$
give endomorphisms of the spinor bundle $\tilde S\otimes L$. Via
pullback and the isomorphism of $S^+\otimes L$ and
$S^-\otimes L$, a 1-form $\alpha$ on $Y$ acts on $\tilde S\otimes L$
as the 2-form $\alpha\wedge dt$ acts on $S^+\otimes L$.

Write $F_A^+=\frac{1}{2}(F_A+*F_A)$, and $F_A=dA$ in coordinates. Using
the fact that $F^-_A$ acts trivially on $S^+\otimes L$ and what was pointed
out in the previous sentence, convince yourself
that the action of $F_A^+$ on $S^+\otimes L$ corresponds exactly to
the action of $\frac{dA}{dt}+*F_A$ on $\tilde S\otimes L$ (in the latter
the Hodge $*$ operator is taken on $Y$).

An analogous argument explains the presence of the term $<e_i\psi,\psi>e^i$.

\noindent QED

The next step in constructing a Seiberg--Witten theory on three-manifolds
is to interpret the above equations (\ref{3eqSW1}), (\ref{3eqSW2}), as
the downward gradient flow of a functional defined on connections and
sections of the spinor bundle of $Y$. With this interpretation, solutions to
the Seiberg--Witten equations will be critical points of the functional
and the associated invariant will admit an interpretation as an
Euler characteristic of a suitable space of connections and sections.

\begin{defin}
Let's introduce a functional, defined on the space of connections on $L$
and sections of $\tilde S\otimes L$ over $Y$.
\beq
C(A,\psi)=\frac{1}{2}\int_Y (A-A_0)\wedge F_A +
\frac{1}{2}\int_Y <\psi,\partial_A\psi>dv.
\label{3funct}
\eeq
\end{defin}

\begin{teor}
The gradient flow of the functional $C$ is given by
\[ (\partial_A\psi, *F_A-<e_i\psi,\psi>e^i). \]
\label{flow}
\end{teor}

We would like to have the functional (\ref{3funct}) defined on connections and
sections, modulo the action of the gauge group. However, (\ref{3funct}) is
invariant under gauge transformations connected to the identity; but, in
the case of maps that belong to other connected components in
\[ \pi_0({\cal G})=H^1(Y\times [0,1];{\bf Z}), \]
an easy computation \cite{KM1} shows that the functional changes according to
\beq
C(A-2i\lambda^{-1}d\lambda,i\lambda\psi)=C(A,\psi)+4\pi^2
<c_1(L)\cup [\lambda],[Y]>,
\label{gautrans}
\eeq
where $[\lambda]$ is the cohomology class of the 1-form
$\frac{-i\lambda^{-1}}{\pi}d\lambda$, representing the connected component
of $\lambda$ in the gauge group.

The functional (\ref{3funct}) represents an analogue of the Chern--Simons
functional in Donaldson's gauge theory.
As in the case of the Chern-Simons form, the gradient flow defines a vector
field on the space $\hat{\cal A}/{\cal G}$ of connections and sections
of $\tilde S\otimes L$ over the 3-manifold $Y$.  We assume to have only
non-degenerate solutions, $\psi\neq 0$.

Solutions of (\ref{3eqSW1}), (\ref{3eqSW2}), are critical points of this
vector field. The Hessian at a critical point can be described as the
determinant of a Fredholm operator that linearizes the equations.

In fact, we have an analogue of the complex (\ref{chcomplex}) given by the
infinitesimal action of the gauge group and the linearization of
the equation:
\beq
0\rightarrow \Lambda^0(Y)\stackrel{G}{\rightarrow} \Lambda^1(Y)\oplus
\Gamma(Y,\tilde S\otimes L)\stackrel{T}{\rightarrow} \Lambda^1(Y)
\oplus \Gamma(Y,\tilde S\otimes L)\rightarrow 0,
\label{3chcomplex}
\eeq
with the map $T$ given by
\[ T(\alpha,\phi)=\left( \begin{array}{c}
                  -*d\alpha +2Im(<e_i\psi_0,\phi>)e^i\\
                  -\partial_{A_0}\phi-i\alpha\cdot\psi_0
                          \end{array} \right). \]

The operator $G^*+T$ is elliptic; up to perturbing the equations we
can ensure that minus the Euler characteristic of the complex $C^*$ is given
just by the first Betti number, which represents therefore the generic
dimension of the moduli space.

When $Ker(G^*+T)=H^1(C^*)=0$, the moduli
space consist generically of finitely many points, and the invariant is
defined as the sum with the orientation given by the determinant of the
Fredholm operator.

This is related to the Hessian of the functional
(\ref{3funct}) by the definition of the index of a critical point, in
an infinite dimensional context, with the spectral flow.

Using the complex (\ref{3chcomplex}) and the Mathai--Quillen
formalism, it is possible to construct a topological Lagrangian for the
three-dimensional theory as well. This is dealt with in \cite{ZWC}. The
upshot is that the invariant can be recovered as a partition function of
the QFT. From this point of view it seems that such an invariant should
represent a Seiberg--Witten analogue of the Casson invariant, whose gauge
theoretic interpretation is explained in \cite{T7}.

It would be interesting to have a combinatorial, or better knot-theoretic,
description of this invariant. It is rumored that Turaev has done something
in this direction, but we don't have, at present, any precise information.

An application of the Seiberg--Witten techniques to knot theory can be
found in \cite{Au}. The estimate (\ref{est}) used in step (4)
of the proof of the Thom conjecture by Kronheimer and Mrowka
is applied to prove that the knot T(5,2) has unknotting number 2.

\subsection{Exercises}
\begin{itemize}

\item Try to give a more detailed proof of theorem \ref{dimred}.

\item Sketch a proof of theorem \ref{flow}: compute the variation
of the functional, $\delta C$, corresponding to an increment $A+\epsilon e^i$,
$\psi +\epsilon\phi$.

\item Complete the proof of point (3) of the Thom conjecture: following
\cite{KM1}, consider solutions $(A_R,\psi_R)$ on the cylinder
$Y\times [-R, R]$; show that the change in the functional $C(A_R(R),
\psi_R(R))-C(A_R(-R),\psi_R(-R))$ is negative and uniformly bounded,
independently of $R$. For the latter property consider gauge transformations
such that $A_R-2i\lambda_R^{-1}d\lambda_R$ and the
first derivatives are uniformly bounded. Show that the functional $C$
changes monotonically along the cylinder, and deduce that for all $N$ there
is a solution on $Y\times [0,1]$ for which the change of $C$ is
bounded by $1/N$. Complete the argument by showing that there is a
translation invariant solution.

\item How can the ``stretching the neck'' argument for the connected
sum theorem work? Sketch a proof, by considering solutions on cylinders
$S^3\times [-R,R]$. Use information on the scalar curvature.

\end{itemize}

\section{Towards a Topological Interpretation of Quantum Co\-homology?}

In the world of symplectic manifolds, there are invariants defined in terms
of the {\em pseudo-holomorphic curves}, i.e. of those 2-submanifolds embedded
via
a map which is holomorphic with respect to a Riemann surface structure on
the domain and an almost-complex structure $J$ tamed by the symplectic form
on the manifold $X$. These are known as Gromov invariants, \cite{McDS}.
A symplectic form $\omega$ is said to tame an almost-complex structure $J$ if
it satisfies $\omega(v, Jv)>0$ for all non zero vectors $v$.

A recent paper by Taubes \cite{T2} has uncovered a deep relation between
these and the Seiberg-Witten invariants.

In our context, we consider the case of a 4-dimensional compact connected
symplectic manifold without boundary.

Given a homology class $A\in H_2(X;{\bf Z})$, let ${\cal H}(A)$ be the
space of $J$-holo\-morphic curves that realize the class $A$.

\begin{lem}
For a generic choice of the almost-complex structure $J$, the space
${\cal H}(A)$ is a smooth manifold of dimension
\[ d_{\cal H}=<-c_1(K)\cup PD(A)+PD(A)^2,[X]>, \]
where $K$ is the canonical line bundle defined by the symplectic
structure and $PD(A)$ is the Poicar\'e dual of the class $A$.
\label{dimhol}
\end{lem}

This result is proven with the same technique we have used in computing the
dimension of the Seiberg--Witten moduli space: in fact the equation
describing the $J$-holomorphic condition linearizes to a Fredholm
operator.

When $d_{\cal H}>0$, given a set $\Omega$ of $\frac{d_{\cal H}}{2}$ distinct
points in $X$, we take ${\cal H}_\Omega$ to be the subspace of ${\cal H}(A)$
whose points are the curves that contain all of the points in $\Omega$.

\begin{teor}
For a generic choice of $J$ and of the points in $X$, ${\cal H}_\Omega$
is a compact zero--dimensional manifold endowed with a canonical
orientation.
\label{holzero}
\end{teor}

Some details can be found in the first chapter of \cite{McDS}.
This allows us to define invariants as follows.

\begin{defin}
We define the Gromov invariant to be a map
\[ \Phi: H_2(X;{\bf Z})\to {\bf Z} \]
which is the sum with orientation of the points in ${\cal H}_\Omega(A)$,
if $d_{\cal H}>0$; the sum of points of ${\cal H}(A)$ with orientation, if
$d_{\cal H}=0$; and zero by definition when $d_{\cal H}<0$.
\label{grominv0}
\end{defin}

The value of $\Phi$ is independent of the choice of a generic set of
points and of the quasi-complex structure $J$.

The result announced in Taubes' paper \cite{T2} is the following.

\begin{teor}
Given a line bundle $L^2$ on a compact symplectic 4-manifold $X$, we have
\[  N_L\equiv \Phi(PD(c_1(L^2))). \]
\label{SWG}
\end{teor}

The strategy of the proof is to extend to other $Spin_c$-structures the
asymptotic technique introduced in \cite{T} to compute the Seiberg-Witten
invariant of a symplectic manifold with respect to the canonical
$Spin_c$-structure.

In this process, the asymptotic method shows that to a solution of the
Seiberg-Witten equations corresponds the zero set of a section of a
holomorphic bundle, which is a pseudo-holomorphic curve.

A different argument, which is briefly sketched in \cite{T2}, provides
a converse construction of a solution to the equations, given a $J$-holomorphic
curve. A combined use of these constructions would complete the proof
of the theorem.

A more detailed argument for the first part of the proof has appeared
in a subsequent paper by Taubes \cite{T8}. The second and third parts are
in preparation.
Note that a detailed proof of the theorem is far from being
trivial: the aforementioned paper \cite{T8} is more
than a hundred pages long.

A very interesting problem arises in connection with this result.

The question is whether the connection between the Seiberg--Witten and the
Gromov invariants may lead to a topological formulation of Quantum
Cohomology.

The quantum cup product is a deformed product in the cohomology ring of
a symplectic manifold. Geometrically it can be thought of as a coarser
notion of intersection of homology classes realized by embedded
submanifolds.

Instead of counting, with the orientation, the number of intersection
points of two cycles in generic position, the counting is made over
the $J$-holomorphic curves that touch the given cycles in generic points.

Note that to define this product there are several technical hypothesis
to introduce, which are explained in \cite{McDS}. In our case, since
we deal with the 4-dimensional case only, the situation is simpler.

More precisely, the quantum cohomology ring of $X$ is
\[ QH^*(X)=H^*(X)\otimes {\bf Z}[q,q^{-1}], \]
where $q$ is a formal variable of degree $2N$, with $N$ the minimal
Chern number of $X$ (see \cite{McDS}). A class $a\in QH^k(X)$ splits
as $a=\sum_i a_iq^i$, with $a_i\in H^{k-2iN}(X)$.

There is a non-degenerate pairing
\[ QH^*(X)\otimes QH^*(X)\to {\bf Z}, \]
\[ <a,b>=\sum_{2(i+j)N=k+l-2n} a_ib_j, \]
with $k=deg(a)$, $l=deg(b)$ in $QH^*(X)$, and $n=dim(X)$.

The quantum cup product
\[ QH^k(X)\otimes QH^m(X)\stackrel{*}{\rightarrow} QH^{l+m}(X) \]
is defined by specifying the values of $<a*b,c>$:
\[ <a*b,c>:=\sum_{i,j,k}\sum_A \Phi_A(\alpha_i,\beta_j,\gamma_k), \]
with $a=\sum a_iq^i$, $b=\sum b_jq^j$, $c=\sum c_kq^k$,
$\alpha_i=PD(a_i)$, $\beta_j=PD(b_j)$, $\gamma_k=PD(c_k)$. $A$ is a
homology class realized by a $J$-holomorphic curve of genus zero, with
$c_1(A)+(i+j+k)N=0$, $c_1(A)$ being the evaluation over $A$ of the
Chern class of the restriction of the tangent bundle over $A$ (the latter
condition is imposed for dimensional reasons).

The coefficient $\Phi_A(\alpha_i,\beta_j,\gamma_k)$ is a more general
version of the Gromov invariants introduced in definition \ref{grominv0}.

Here, instead of imposing $d_{\cal H}/2$ points in ${\cal H}_\Omega$,
we want a positive dimensional manifold. Hence
$\Omega=\{ p_1,\ldots,p_r \}$, $r< d_{\cal H}/2$.

Now we would like to map this manifold to $X^r$, so that it gives rise
to a pseudo-cycle (see \cite{McDS}). Thus, we consider an evaluation map
\[ \epsilon_r : {\cal H}_\Omega \to X^r \]
whose image is the space
\[ \{ (u(p_1),\ldots, u(p_r)) \mid u:\Sigma\to X, [u(\Sigma)]=A \}, \]
where $u$ is a $J$-holomorphic parametrisation of the curve.

It is clear, however, that we need to use parametrised
curves. Hence the space ${\cal H}_\Omega$ has to be intended rather
as the moduli space of parametrised curves that touch the points in
$\Omega$, modulo reparametrisations, i.e. modulo automorphisms
of the Riemann surface $\Sigma$.

The topology of this moduli space is more complicated than the
previous case: in fact there are non-trivial problems related to the
compactification. The theory works sufficiently well only in the case of
genus zero curves (see \cite{McDS}, chapter 5).

Once a pseudo-cycle of a certain dimension (which is computed in
\cite{McDS}, chapter 7) is defined in $X^r$, the Gromov invariant
is obtained by intersecting this cycle in $X^r$ (with the usual
intersection product) with a number of homology classes of $X$,
so to reach the complementary dimension.

Although the definition of this invariant is different from the one
used in \cite{T2}, it would be very interesting to describe it in terms of
Seiberg--Witten theory. In fact that could be a step towards a generalisation
of the quantum product to non-symplectic manifolds, and could lead to
a topological, or gauge theoretic, interpretation of quantum cohomology.

One of the highly non-trivial results in quantum cohomology is the fact that
the quantum cup product is associative. It would be really neat to have a
different proof based on these gauge theoretic techniques. But this is just
pure speculation.

\section{Seiberg--Witten and Donaldson's Theory}

Here we would like to mention, without giving any detail,
the conjecture about the relationship between Seiberg--Witten and Donaldson
theory.

In a fundamental paper \cite{KM2}, Kronheimer and Mrowka gave a description
of a relation that constrains the values of the Donaldson invariants for
a manifold of {\em finite type}. The finite type assumption is a technical
hypothesis on the behaviour of the polynomial invariants, which is
satisfied by all known examples of simply connected 4-manifolds.

The result states that the polynomial invariants \cite{D} of a simply
connected manifold of finite type are recovered form the expression
\[ q=\exp(\frac{Q}{2})\sum_{k=1}^p a_k e^{x_k}, \]
where $q$ is the Donaldson invariant
\[ q=\sum_d \frac{q_d}{d!}, \]
\[ q_d(\alpha)=<\omega(\alpha),[M_d]>. \]
Here $M_d$ is the moduli space of Donaldson's theory; and the cohomology
class $\omega(\alpha)$ is constructed in a way (\cite{D}, \cite{DK})
similar to what we've been discussing in the chapter on QFT.
On the right hand side we have the intersection form $Q$, and a finite
number of classes, the {\em basic classes},
\[ x_k\in H^2(X;{\bf Z}) \]
such that the mod 2 reduction of each $x_k$ is the class $w_2(X)$,
and non-zero rational coefficients $a_k$.

The conjecture formulated by Witten in \cite{W} is the following.

\begin{conj}
In the above expression, the basic classes $x\in H^2(X;{\bf Z})$
are exactly those that satisfy
\[ x^2=c_1(L)^2=\frac{2\chi +3\sigma}{4}, \]
i.e. those $Spin_c$-structures that give rise to a zero dimensional
Seiberg--Witten moduli space.
The corresponding coefficient $a_x$ is exactly the Seiberg--Witten
invariant $N_L$.
\label{SWDon}
\end{conj}

A sketch of the Physics underlying this conjecture can be found in
\cite{W} or in \cite{W3}. More detailed references are \cite{WS1},
\cite{WS2}.

{}From the conjecture \ref{SWDon} it seems that the Seiberg--Witten
invariants should contain more information than the
Donaldson invariants. In fact all the Donaldson polynomials can be
recovered from the knowledge of the Seiberg--Witten invariants
associated to zero dimensional moduli spaces. Thus, in principle, the
Seiberg--Witten invariants associated to positive dimensional moduli spaces
might give more information.

A related conjecture, which we mentioned already, is the following.

\begin{conj}
For a manifold $X$ of simple type the only non-trivial Seiberg--Witten
invariants correspond to a choice of $L$ such that
$d=dim(M_L)=0$.
\end{conj}

A first important step towards the proof of the Witten conjecture \ref{SWDon}
has been made recently by Pidstrigach and Tyurin \cite{PiT2}.

\section{Non--abelian Monopoles}

A generalisation of the monopole equations to the case of
a non-abelian structure group has been investigated \cite{LabM2}.

Let us assume, for simplicity, that the manifold $X$ is $Spin$ in order
to avoid the problem of definiteness of the spinor bundles $S^\pm$.

Let $G$ be a compact connected simple Lie group. Instead of the line bundle
$L$, consider a principal $G$-bundle $P$ and an associated vector bundle $E$
obtained via a representation of $G$.

The indeterminates will be $(A,\psi)$ with $A$ a $G$-connection on $E$ and
$\psi\in\Gamma(X,S^+\otimes E)$.

The gauge group ${\cal G}$ of $E$ acts via the representation as
\[ \lambda(A,\psi)=(\lambda^{-1}A\lambda-2\lambda^{-1}d\lambda,
\lambda\psi). \]

The non-abelian monopole equations are
\[ \begin{array}{c}
       D_A\psi=0 \\
       F_A^+ -\frac{1}{4}<e_ie_j\psi,\psi>e^i\wedge e^j =0,
   \end{array} \]
where $D_A$ is the Dirac operator twisted with the $G$-connection $A$,
and the 2-form in the second equation takes values in the gauge Lie
algebra bundle $\Lambda^{2+}({\bf g}_E)$. In fact the inner product of
sections in the second equation can be written as
\[ <e_ie_j\tilde\psi\otimes s,\tilde\psi\otimes s>_{S^+} \in
\Lambda^0({\bf g}_E), \]
where $\psi=\tilde\psi\otimes s$.

The infinitesimal action of the Lie group at a point $(A,\psi)$,
is given by
\[ u\stackrel{G}{\mapsto} (-d_A u,iu\psi), \]
with $\lambda=\exp(iu)$, $u\in\Lambda^0({\bf g}_E)$,
and $d_A=d+i[A,\cdot ]$.

The linearization of the non-abelian monopole equations fits
into a chain complex with the infinitesimal action of the gauge group
\cite{LabM2}:
\[ 0\rightarrow\Lambda^0({\bf g}_E)\stackrel{G}{\rightarrow}\Lambda^1
({\bf g}_E)\oplus\Gamma(X,S^+\otimes E)\stackrel{T}{\rightarrow}\Lambda^{2+}
({\bf g}_E)\oplus\Gamma(X,S^-\otimes E)\rightarrow 0.\]

The operator $T+G^*$ is elliptic and, up to zero order terms, has index
determined by the index of the twisted Dirac operator $D_A$ and by
the complex
\[ 0\rightarrow \Lambda^0({\bf g}_E){\rightarrow}\Lambda^1({\bf g}_E)
{\rightarrow}\Lambda^{2+}({\bf g}_E)\rightarrow 0. \]

For the latter see \cite{AHS}. In the case of $G=SU(N)$ the dimension
of the moduli space turns out to be
\[ dim(M)=(4N-2)c_2(E)-\frac{N^2-1}{2}(\chi +\sigma)-\frac{\delta}{4}
\sigma, \]
where $\delta$ is the dimension of the representation.

The orientability proof also carries over following \cite{D2}, while
the compactness argument fails for the non-abelian case.

In fact, even more can be said \cite{LabM}. In the case of $SU(N)$
consider the singular solutions corresponding to $\psi\equiv 0$. The
second equation becomes the well known anti-self-dual equation of
Donaldson theory.

Hence the whole Donaldson moduli space is recovered as a singular
submanifold of the non-abelian monopole moduli space.

The result is interesting as a step towards a better understanding of the
relation between Donalson and Seiberg--Witten invariants.
In fact non-abelian monopoles are a major ingredient in the approach
to the Witten conjecture \ref{SWDon} by Pidstrigach and
Tyurin \cite{PiT2}.

Another interesting result, that seems strictly related to the aforementioned,
is obtained by purely physical arguments in \cite{BL}. The authors
construct a unitary transformation between the QFT operators that give
the Floer homology \cite{Fl} of Donaldson theory and those obtained from
the non-abelian monopole equation. They conjecture that a better understanding
of this relation could provide new insight on the Witten conjecture
\ref{SWDon}.

\end{document}